\documentclass[11pt]{article}
\usepackage[colorlinks=true, linkcolor=black, citecolor=black]{hyperref}
\usepackage[margin=1cm]{caption}
\usepackage{amssymb}
\usepackage{amsthm}
\usepackage{float,url}
\usepackage{tikz}
\usepackage[left=1.00in, right=1.00in, top=1.00in, bottom=1.00in]{geometry}
\usepackage{mathdots}
\def\chil{\chi_{\ell}}
\newcommand\floor[1]{\lfloor#1\rfloor}
\newcommand\Floor[1]{\left\lfloor#1\right\rfloor}

\newcommand\Ceil[1]{\left\lceil#1\right\rceil}
\author{Daniel W. Cranston\thanks{
Department of Mathematics and Applied
Mathematics, Virginia Commonwealth University, Richmond, VA, 23284. email:
\texttt{dcranston@vcu.edu}} \and Bobby Jaeger\thanks{
Department of Mathematics and Applied
Mathematics, Virginia Commonwealth University, Richmond, VA, 23284. email:
\texttt{jaegerrj@vcu.edu}}
}

\def\ch{\textrm{ch}}
\newtheorem*{mainthm}{Main Theorem}
\newtheorem*{basic}{Basic Reducibility Lemma}
\newtheorem*{reducibility}{Main Reducibility Lemma}
\newtheorem*{concavitylem}{Concavity Lemma}
\newtheorem{cor}{Corollary}
\newtheorem{conj}{Conjecture}
\newtheorem{theorem}{Theorem}[section]
\newtheorem{conjecture}[theorem]{Conjecture}

\def\VertSize{0.1}

\begin{document}

\title{List-coloring the Squares of Planar Graphs\\ without 4-Cycles and 5-Cycles}
\maketitle
\begin{abstract}
Let $G$ be a planar graph without 4-cycles and 5-cycles and with maximum degree
$\Delta\ge 32$.  We prove that $\chil(G^2)\le \Delta+3$. For arbitrarily large
maximum degree $\Delta$, there exist planar graphs $G_{\Delta}$ of girth 6 with
$\chi(G_{\Delta}^2)=\Delta+2$.  Thus, our bound is within 1 of being optimal.  Further,
our bound comes from coloring greedily in a good order, so the bound
immediately extends to \emph{online} list-coloring.  In addition, we prove
bounds for $L(p,q)$-labeling. Specifically, $\lambda_{2,1}(G)\le \Delta+8$ and,
more generally, $\lambda_{p,q}(G)\le (2q-1)\Delta+6p-2q-2$, for positive
integers $p$ and $q$ with $p\ge q$.  Again, these bounds come from a
greedy coloring, so they immediately extend to the list-coloring and online
list-coloring variants of this problem.
\end{abstract}

\section{Introduction}

The \emph{square} $G^2$ of a graph $G$ is formed from $G$ by adding an edge between each
pair of vertices at distance two in $G$.
In 1977, Wegner \cite{wegner} posed the following conjecture, which has attracted great
interest, and led to a remarkable number of results.
(Most of our terminology and notation is standard.  When it is not, we define
terms where they are first used.  For reference, we also collect some key
definitions in the \hyperref[appendix]{Appendix}.)

\begin{conjecture}[Wegner \cite{wegner}]
If $G$ is a planar graph with maximum degree $\Delta$, then
\[ \chi(G^2) \le 
\left\{
\begin{array}{ll}
7 & \mbox{if } \Delta=3;\\
\Delta+5 & \mbox{if } 4\le\Delta\le7;\\
\floor{\frac{3\Delta}{2}} + 1 & \mbox{if } \Delta \ge 8.
\end{array}
\right.
\]
\end{conjecture}

\begin{figure}[h]
\centering
\begin{tikzpicture}[scale=1.2]
\coordinate (A) at (0, 3);
\coordinate (B) at (-1.73205, 0);
\coordinate (C) at (1.73205, 0);

\draw (A) to[out = -150, in = 90] (B) to[out = 30, in = -90] (A) to[out = -30, in = 90] (C) to[out = -150, in = -30] (B) to[out = 30, in = 150] (C) to[out = 150, in = -90] (A);
\draw (B) to[out = -90, in = 180] (0, -1.3) to[out = 0, in = -90] (C);

\draw[fill] (A) circle [radius = \VertSize];
\draw[fill] (B) circle [radius = \VertSize];
\draw[fill] (C) circle [radius = \VertSize];

\draw[fill] (-1.275, 1.8) circle [radius = \VertSize];
\draw[fill] (1.275, 1.8) circle [radius = \VertSize];
\draw[fill] (-0.45, 1.2) circle [radius = \VertSize];
\draw[fill] (0.45, 1.2) circle [radius = \VertSize];
\draw[fill] (0, 0.525) circle [radius = \VertSize];
\draw[fill] (0, -0.525) circle [radius = \VertSize];

\node at (0,3.23) {$u$};
\node[left] at (B) {$v\,$};
\node[right] at (C) {$\,w~$};

\node at (-0.86603, 1.65) {$ \ddots $};
\node at (0.86603, 1.65) {$ \iddots $};
\node at (0, 0.15) {$ \vdots $};

\node[rotate=50,yscale=2.0] at (-0.7, 1.75) {\}};
\node at (-0.55, 2.1) {$\left\lceil\frac{\Delta}{2}\right\rceil$};
\node[rotate=-50,yscale=2.0] at (0.7, 1.75) {\{};
\node at (0.55, 2.1) {$\left\lfloor\frac{\Delta}{2}\right\rfloor$};
\node[yscale=2.0] at (0.2, 0) {\}};
\node at (0.8, 0) {$\left\lfloor\frac{\Delta}{2}\right\rfloor - 1$};
\end{tikzpicture}
\caption{Wegner's construction for $\Delta \ge 8$.}
\label{fig:fat}
\end{figure}
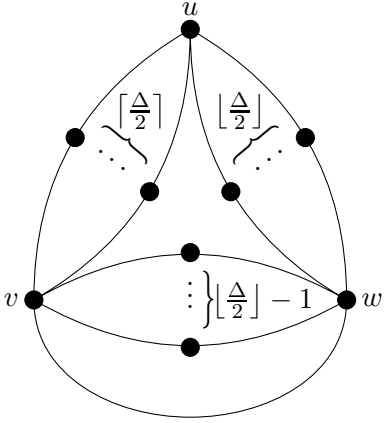

Wegner also gave constructions showing that this conjecture is sharp if true.
In particular, his sharpness example for $\Delta \ge 8$ is shown in
Figure~\ref{fig:fat}. 
Although the conjecture remains open in general, Havet et al.\  
\cite{havet} showed that the conjectured upper bound
holds asymptotically, i.e., $\chi(G^2) \le \frac{3}{2}\Delta + o(\Delta)$. 
A more thorough history of Wegner's conjecture appears in the introductions of
\cite{dvorak} and \cite{havet}.

For every graph $G$, we have the lower bound $\chi(G^2)\ge \Delta+1$.
If we seek to prove an upper bound closer to this trivial lower bound, we
clearly must forbid the configuration of Figure~\ref{fig:fat}.  
Forbidding 3-cycles alone does not
really help, since now subdiving the edge $vw$ yields a graph $G$ with no
3-cycles and such that $G^2$ still has clique number $\Floor{\frac32\Delta}$.
So we make the obvious choice and forbid 4-cycles, as well as perhaps cycles of
other lengths.
This line of inquiry has an intriguing history, much of which was motivated by
the following conjecture of Wang and Lih~\cite{wanglih}.

\begin{conj}[Wang and Lih~\cite{wanglih}]
For every integer $g$ at least 5, there exists some integer $\Delta_g$ such that
every planar graph $G$ with girth at least $g$ and maximum degree at least
$\Delta_g$ satisfies $\chi(G^2)=\Delta+1$.
\label{conj1}
\end{conj}

The conjecture was proved by Borodin et al.~\cite{girthseven} for $g\ge 7$
and disproved for $g\in \{5,6\}$ in the same paper.  However, Dvo\v{r}\'{a}k et
al.~\cite{dvorak} complemented these results with the following theorem.

\begin{theorem}[\cite{dvorak}]
If $G$ is a planar graph with girth at least 6 and $\Delta\ge 8821$, then
$\chi(G^2)\le \Delta+2$.
\end{theorem}

(Soon after, Borodin et al.~\cite{girthsix2} weakened the hypothesis to
$\Delta\ge 18$.) 
In the same paper, Dvorak et al.~posed the following conjecture.

\begin{conj}
There exists some constant $M$ such that every planar graph $G$ with girth 5 and
maximum degree at least $M$ satisfies $\Delta(G^2)\le \Delta+2$.
\label{conj2}
\end{conj}

If true, Conjecture~\ref{conj2} would be a very nice result.  
Zhu et al.~\cite{zhu-etal} went in a slightly different direction.  They
considered planar graphs with no 4-cycles and no 5-cycles (although 3-cycles are allowed).
Among other results, they showed that if $\Delta\ge 9$, then $\chi(G^2)\le
\Delta+5$.  In fact, this bound follows from a more general result on
$L(p,q)$-labeling, which we will discuss soon.

Our main result is the following theorem.

\begin{mainthm}
Let $G$ be a planar graph with maximum degree $\Delta$ that contains no
4-cycles and no 5-cycles.  If $\Delta \ge 32$, then there exists an ordering
$v_1, \ldots, v_n$ of $V(G)$ such that each $v_i$ has at most 3 neighbors
in $G$ that appear earlier in the ordering and at most $\Delta+2$ neighbors in
$G^2$ that appear earlier in the ordering.
\label{main}
\end{mainthm}

This theorem is optimal in the following sense.
We cannot reduce the bound of ``at most 3 neighbors in $G$''
to ``at most 2''.  To see this, it suffices to construct planar graphs with
arbitrarily large maximum degree, no 4-cycles and no 5-cycles, and minimum
degree 3.  We do so as follows.  

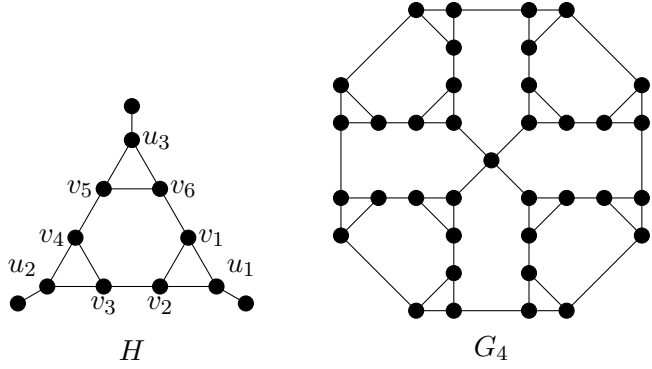
\begin{figure}
\centering
\begin{tikzpicture}
\coordinate (A) at (0.75, 0.0);
\coordinate (B) at (0.375, -0.64951);
\coordinate (C) at (-0.37499, -0.64951);
\coordinate (D) at (-0.75, 0.0);
\coordinate (E) at (-0.375, 0.64951);
\coordinate (F) at (0.375, 0.64951);
\coordinate (G) at (0.0, 1.29903);
\coordinate (H) at (1.125, -0.64951);
\coordinate (I) at (-1.12499, -0.64951);
\coordinate (J) at (0.0, 1.75);
\coordinate (K) at (1.51554, -0.87499);
\coordinate (L) at (-1.51554, -0.875);

\draw (A)--(B)--(C)--(D)--(E)--(F)--(A) (A)--(H)--(B) (C)--(I)--(D) (E)--(G)--(F) (G)--(J) (H)--(K) (I)--(L);

\draw[fill] (A) circle [radius = \VertSize];
\draw[fill] (B) circle [radius = \VertSize];
\draw[fill] (C) circle [radius = \VertSize];
\draw[fill] (D) circle [radius = \VertSize];
\draw[fill] (E) circle [radius = \VertSize];
\draw[fill] (F) circle [radius = \VertSize];
\draw[fill] (G) circle [radius = \VertSize];
\draw[fill] (H) circle [radius = \VertSize];
\draw[fill] (I) circle [radius = \VertSize];
\draw[fill] (J) circle [radius = \VertSize];
\draw[fill] (K) circle [radius = \VertSize];
\draw[fill] (L) circle [radius = \VertSize];

\node[right] at (A) {$ v_1 $};
\node[below] at (B) {$ v_2 $};
\node[below] at (C) {$ v_3 $};
\node[left] at (D) {$ v_4 $};
\node[left] at (E) {$ v_5 $};
\node[right] at (F) {$ v_6 $};
\node[above right] at (H) {$ u_1 $};
\node[above left] at (I) {$ u_2 $};
\node[right] at (G) {$ u_3 $};

\node at (0, -1.5) {$ H $};
\end{tikzpicture}
~~~~~
\begin{tikzpicture}
\coordinate (A1) at (1, 0.5);
\coordinate (B1) at (1.5, 0.5);
\coordinate (C1) at (2, 1);
\coordinate (D1) at (1, 2);
\coordinate (E1) at (0.5, 1.5);
\coordinate (F1) at (0.5, 1);
\coordinate (G1) at (0.5, 0.5);
\coordinate (H1) at (2, 0.5);
\coordinate (I1) at (0.5, 2);

\draw (A1)--(B1)--(C1)--(D1)--(E1)--(F1)--(A1) (A1)--(G1)--(F1) (B1)--(H1)--(C1) (D1)--(I1)--(E1);

\draw[fill] (A1) circle [radius = \VertSize];
\draw[fill] (B1) circle [radius = \VertSize];
\draw[fill] (C1) circle [radius = \VertSize];
\draw[fill] (D1) circle [radius = \VertSize];
\draw[fill] (E1) circle [radius = \VertSize];
\draw[fill] (F1) circle [radius = \VertSize];
\draw[fill] (G1) circle [radius = \VertSize];
\draw[fill] (H1) circle [radius = \VertSize];
\draw[fill] (I1) circle [radius = \VertSize];

\coordinate (A2) at (-1, 0.5);
\coordinate (B2) at (-1.5, 0.5);
\coordinate (C2) at (-2, 1);
\coordinate (D2) at (-1, 2);
\coordinate (E2) at (-0.5, 1.5);
\coordinate (F2) at (-0.5, 1);
\coordinate (G2) at (-0.5, 0.5);
\coordinate (H2) at (-2, 0.5);
\coordinate (I2) at (-0.5, 2);

\draw (A2)--(B2)--(C2)--(D2)--(E2)--(F2)--(A2) (A2)--(G2)--(F2) (B2)--(H2)--(C2) (D2)--(I2)--(E2);

\draw[fill] (A2) circle [radius = \VertSize];
\draw[fill] (B2) circle [radius = \VertSize];
\draw[fill] (C2) circle [radius = \VertSize];
\draw[fill] (D2) circle [radius = \VertSize];
\draw[fill] (E2) circle [radius = \VertSize];
\draw[fill] (F2) circle [radius = \VertSize];
\draw[fill] (G2) circle [radius = \VertSize];
\draw[fill] (H2) circle [radius = \VertSize];
\draw[fill] (I2) circle [radius = \VertSize];

\coordinate (A3) at (-1, -0.5);
\coordinate (B3) at (-1.5, -0.5);
\coordinate (C3) at (-2, -1);
\coordinate (D3) at (-1, -2);
\coordinate (E3) at (-0.5, -1.5);
\coordinate (F3) at (-0.5, -1);
\coordinate (G3) at (-0.5, -0.5);
\coordinate (H3) at (-2, -0.5);
\coordinate (I3) at (-0.5, -2);

\draw (A3)--(B3)--(C3)--(D3)--(E3)--(F3)--(A3) (A3)--(G3)--(F3) (B3)--(H3)--(C3) (D3)--(I3)--(E3);

\draw[fill] (A3) circle [radius = \VertSize];
\draw[fill] (B3) circle [radius = \VertSize];
\draw[fill] (C3) circle [radius = \VertSize];
\draw[fill] (D3) circle [radius = \VertSize];
\draw[fill] (E3) circle [radius = \VertSize];
\draw[fill] (F3) circle [radius = \VertSize];
\draw[fill] (G3) circle [radius = \VertSize];
\draw[fill] (H3) circle [radius = \VertSize];
\draw[fill] (I3) circle [radius = \VertSize];

\coordinate (A4) at (1, -0.5);
\coordinate (B4) at (1.5, -0.5);
\coordinate (C4) at (2, -1);
\coordinate (D4) at (1, -2);
\coordinate (E4) at (0.5, -1.5);
\coordinate (F4) at (0.5, -1);
\coordinate (G4) at (0.5, -0.5);
\coordinate (H4) at (2, -0.5);
\coordinate (I4) at (0.5, -2);

\draw (A4)--(B4)--(C4)--(D4)--(E4)--(F4)--(A4) (A4)--(G4)--(F4) (B4)--(H4)--(C4) (D4)--(I4)--(E4);

\draw[fill] (A4) circle [radius = \VertSize];
\draw[fill] (B4) circle [radius = \VertSize];
\draw[fill] (C4) circle [radius = \VertSize];
\draw[fill] (D4) circle [radius = \VertSize];
\draw[fill] (E4) circle [radius = \VertSize];
\draw[fill] (F4) circle [radius = \VertSize];
\draw[fill] (G4) circle [radius = \VertSize];
\draw[fill] (H4) circle [radius = \VertSize];
\draw[fill] (I4) circle [radius = \VertSize];

\coordinate (X) at (0, 0);

\draw (X)--(G1) (X)--(G2) (X)--(G3) (X)--(G4);
\draw (I1)--(I2) (H2)--(H3) (I3)--(I4) (H4)--(H1);

\draw[fill] (X) circle [radius = \VertSize];

\node at (0, -2.5) {$ G_4 $};
\end{tikzpicture}
\caption{The gadget $H$ (on the left) and $G_4$.}
\end{figure}

Form gadget $H$ from a 6-cycle $v_1\ldots v_6$ by adding vertices $u_1, u_2, u_3$ with
$u_1$ adjacent to $v_1$ and $v_2$; $u_2$ adjacent to $v_3$ and $v_4$; and $u_3$
adjacent to $v_5$ and $v_6$.  Finally, add a pendant edge incident to each
$u_i$.  
To form graph $G_k$,
begin with a cycle $C_k$ and add a dominating vertex.  Now replace,
successively,
each 3-vertex $x$ of the resulting graph with a copy of $H$, joining each
neighbor of $x$ to $H$ using its three pendant edges.  Clearly the resulting
graph has minimum degree 3.  Each cycle within a copy of $H$ has length 3 or at
least 6, and each cycle through more than one copy of $H$ has length at least 9.
Thus, for any ordering $\sigma$ of the vertices of $G_k$, the final vertex will
have at least 3 neighbors earlier in $\sigma$.

To put this theorem in context, we note that this approach of coloring greedily
in a good ordering was used implicitly by van
den Heuvel and McGuinness~\cite{vandenheuvel} in their proof that every planar graph
with $\Delta$ large enough satisfies $\chi(G^2)\le 2\Delta+25$.  The method was
made explicit by Agnarsson and Halld\'{o}rsson~\cite{agnarsson} and Borodin et
al.~\cite{BBGH1,BBGH2} who (independently) improved this result to $\chi(G^2)\le
\Ceil{\frac95\Delta}+1$ for $\Delta$ sufficiently large.  Both groups showed
that this bound is the best possible with this technique, by constructing
planar graphs $G_k$ of arbitrarily high maximum degree $k$ such that $G^2_k$
has minimum degree $\Ceil{\frac95\Delta}$.  This approach has also been used in
some results on $L(p,q)$-labeling.

Our interest in our Main Theorem is due primarily to the following two corollaries.

\begin{cor}
If $G$ is a planar graph with $\Delta\ge 32$ and neither 4-cycles nor 5-cycles,
then $\chil(G^2)\le \Delta+3$.  In fact, this bound holds also for paintability: 
$\chi_p(G^2)\le \Delta+3$.
\end{cor}

The bound on $\chil(G^2)$ comes directly from the Main Theorem, by coloring
greedily in the prescribed ordering.  Since each vertex $v$ has at most
$\Delta+2$ earlier neighbors, some color remains for use on $v$.  For
paintability, the same argument works: on each round, Painter greedily forms a
maximal stable set, by adding vertices in the prescribed order.  As we noted
above, there exist graphs $G_{\Delta}$ with arbitrarily large maximum degree
$\Delta$ for which $\chi(G^2_{\Delta})=\Delta+2$.  (For completeness, we include
in the appendix a construction proving this, due to Dvo\v{r}\'{a}k et
al.~\cite{dvorak}.) Hence, these bounds are within 1 of being best possible.

An $L(p,q)$-labeling is an assignment $f$
of nonnegative
integers to the vertices such that all adjacent vertices $u$ and $v$ satisfy
$|f(u)-f(v)|\ge p$ and vertices $u$ and $v$ at distance two satisfy
$|f(u)-f(v)|\ge q$.  The \emph{$L(p,q)$-labeling number} $\lambda_{p,q}(G)$ is
the minimum value of the largest label $k$ taken over all $L(p,q)$-labelings.
For planar graphs with no 4-cyclces, no 5-cyclces, and $\Delta$ sufficiently
large, Zhu et al.\cite{zhu-etal}~proved that $\lambda_{p,q}\le
(2q-1)\Delta+6p+2q-4$.  In particular, for $\Delta\ge11$, they proved
$\lambda_{2,1}\le \Delta+10$.  In the following corollary, we improve this bound
for $\Delta\ge 32$.
\begin{cor}
If $G$ is a planar graph with $\Delta\ge 32$ and neither 4-cycles nor 5-cycles,
then $\lambda_{p,q}(G)\le (2q-1)\Delta+6p-2q-2$.  In particular,
$\lambda_{2,1}(G)\le \Delta+8$.
\end{cor}

As above, these bounds come from coloring greedily in the prescribed order.
Consider a vertex $v_i$.  Each of its at most 3 earlier neighbors forbid at most
$(2p-1)$ labels; each of its other at most $(\Delta+2-3)$ earlier neighbors in
$G^2$ forbid at most $(2q-1)$ labels.  Since the smallest allowable label is 0,
we get $\lambda_{p,q}(G)\le (2q-1)\Delta+6p-2q-2$.
Note that, also by greedily coloring, 
the bounds generalize immediately to \emph{online list} $L(p,q)$-labeling.

\subsection{Reducibility}
To avoid some technical difficulties (caused by deleting a
vertex and reducing the maximum degree of $G$) we prove the following
theorem, which immediately implies our Main Theorem.

\begin{theorem}
If $G$ is a planar graph with maximum degree $\Delta$ that contains no
4-cycles and no 5-cycles, then there exists an ordering $v_1, \ldots, v_n$ of
$V(G)$ such that each $v_i$ has at most 3 neighbors in $G$ that appear earlier
in the ordering and at most $\max(\Delta,32)+2$ neighbors in $G^2$ that appear
earlier in the ordering.
\label{helper}
\end{theorem}

In what follows, we prove some structural properties of a minimal counterexample
to our theorem. 
Henceforth, let $G$ denote such a minimal counterexample. More
precisely, let $G$ be a planar graph with no 4-cycles and no 5-cycles and
such that no ordering $v_1,\ldots, v_n$ of $V(G)$ has every vertex $v_i$ with
both at most 3 neighbors in $G$ earlier in the ordering and
at most $\max(\Delta,32)+2$ neighbors in $G^2$ earlier in the ordering.  
Moreover, every proper subgraph of $G$ has such an ordering.
Let $N^2(u)$ denote the set of neighbors of $u$ in $G^2$.
Let $D=\max(32,\Delta)$.
We call the ordering guaranteed by the Main Theorem a \emph{good ordering} for
$G$.

\begin{basic}
A minimal $G$ has no vertex $u$ such that $d(u)\le 3$ and $|N^2(u)|\le D+2$ and
$(G-u)^2=G^2-u$.
In particular, (i) $\delta(G)\ge 2$ and (ii) for every 2-vertex $u$
on a 3-cycle $uv_1v_2$ we have $d(v_1)+d(v_2)\ge D+5$.
\label{basicreducibility}
\end{basic}

\begin{proof}
If $u$ is such a vertex, then a good ordering for $G-u$ extends to a good
ordering for $G$ by appending $u$ to the order.  Further, we have $d(u)\le 2$
and $|N^2(u)|\le  D+2$ if either $u$ is (i) a 1-vertex or
(ii) $u$ is a 2-vertex on a 3-cycle $uv_1v_2$ with $d(v_1)+d(v_2)\le D+4$.
\end{proof}

\begin{figure}[h]
\centering
\begin{tikzpicture}
\coordinate (A) at (0, 0);
\coordinate (B) at (1, 0);

\draw (A)--(B);

\draw[fill] (A) circle [radius=\VertSize];
\draw[fill=white] (B) +(-\VertSize, \VertSize) rectangle +(\VertSize, -\VertSize) (B);

\node[left] at (A) {$ u $};
\node[right] at (B) {$ v $};

\node at (0.5, -0.5) {(i)};
\end{tikzpicture}
~~~~~~~~~~~~~~~~~~~~
\begin{tikzpicture}
\coordinate (A) at (0, 1);
\coordinate (B) at (-0.75, 0);
\coordinate (C) at (0.75, 0);

\draw (A)--(B)--(C)--(A);

\draw[fill] (A) circle [radius=\VertSize];
\draw[fill=white] (B) +(-\VertSize, \VertSize) rectangle +(\VertSize, -\VertSize) (B);
\draw[fill=white] (C) +(-\VertSize, \VertSize) rectangle +(\VertSize, -\VertSize) (C);

\node[left] at (A) {$ u $};
\node[left] at (B) {$ v_1 $};
\node[right] at (C) {$ v_2 $};

\node at (0, -0.5) {(ii)};
\end{tikzpicture}
\caption{(i) a 1-vertex is reducible.  (ii) a 2-vertex on a 3-cycle $uv_1v_2$ is
reducible if $d(v_1)+d(v_2)\le  D+4$.}
\label{fig:basic}
\end{figure}
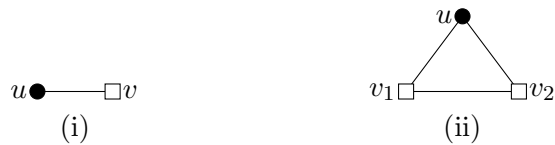

This lemma is illustrated in Figure~\ref{fig:basic}. 
Note that here and throughout the paper,
a vertex that is drawn as a filled circle has all of its incident edges drawn,
while a vertex that is drawn as an empty box may have other incident edges that
are not shown.

We can extend the idea behind the 
\hyperref[basicreducibility]{Basic Reducibility Lemma} to give another,
stronger reducibility lemma.

\begin{reducibility}
\label{mainreducibility}
A minimal $G$ has no
sequence $S = \{w_1, \ldots, w_k\} $ of distinct vertices in $ V(G) $ such that $ E(G[S])
\ne \emptyset$, and also $|N(w_i)\setminus \{w_{i+1},\ldots,w_k\}|\le 3$ and
$|N^2(w_i) \setminus \{w_{i + 1}, \ldots, w_k\}| \le  D + 2$ for every $1 \le i \le k$.
\end{reducibility}

\begin{proof}
Suppose, to the contrary, that such a sequence $S$ exists.
Choose some $e\in E(G[S])$.
Since $G-e$ is a proper subgraph, it has some good ordering $\sigma'$.  To extend
$\sigma'$ to $G$, we delete all elements of $S$ and append them in order; call
this new ordering $\sigma$.  Note that all edges of $G^2$ that are absent from
$(G-e)^2$ are incident with some vertex of $S$.  So $\sigma$ is certainly good
for each vertex of $V(G)\setminus S$.  By hypothesis, it is also good for
each vertex of $S$.
\end{proof}

Whenever we invoke this lemma, we will list the sequence $S$ in the
appropriate order. While this result holds in general, we will typically
use it when $k=2$ or $k=3$.  The case $k=2$ gives the following useful intuition
for the proof: For each edge $u_1u_2$ in $G$, at least one $u_i$ is either
a $5^+$-vertex or has $|N^2(u_i)|\ge  D+3$ (or possibly they are both
4-vertices).  Thus, 
when we do our discharging analysis later,
each edge with an endpoint that needs charge (this will be a vertex of low
degree) has some charge ``nearby'', since it has a nearby vertex of large
degree.  The work of the proof is formalizing this intuition.
\smallskip

To conclude this section, we prove a Concavity Lemma.  Essentially, this lemma
implies that if $|N^2(u)|$ is fixed, then vertex $u$ receives the least charge
when it has one high degree neighbor and all other neighbors have degree as
small as possible (subject to the constraint on $|N^2(u)|$).

\begin{concavitylem}
Let $f(x) = 1 - \frac{4}{x}$, considered on some interval $[a, \infty)$ where
$a > 0$. If $x_1, \ldots, x_n$ are to be chosen in $[a, \infty)$ such that
$\sum_{i = 1}^n x_i = C$ for some constant $C$, then the minimum value of
$\sum_{i=1}^n f(x_i)$ is achieved when $x_1 = \ldots = x_{n - 1} = a$ and $x_n
= C - a(n - 1)$.
\label{concavity}
\end{concavitylem}

\begin{proof}
It suffices to show that $ f(x_1) + f(x_2) \ge f(a) + f(x_1 + x_2 - a) $ for all $ x_1, x_2 \in [a, \infty) $, since we can then proceed by induction on the number of $ x_i $ that are not equal to $ a $.

Assume without loss of generality that $ x_1 \le x_2 $, and let $ t = x_1 - a $. Since $ f $ is concave, its derivative is decreasing, and can be bounded at a point by left and right secants there, giving:
\[ \frac{f(x_2 + t) - f(x_2)}{t} \le f'(x_2) \le f'(x_1) \le \frac{f(x_1) - f(x_1 - t)}{t}. \]
Clearing denominators and rearranging terms gives $ f(x_2 + t) + f(x_1 - t) \le f(x_1) + f(x_2) $. But this is equivalent to $ f(x_1 + x_2 - a) + f(a) \le f(x_1) + f(x_2) $, as was desired.
\end{proof}

\section{Proof of the Main Theorem via Discharging}

Our proof of the Main Theorem is by the discharging method, which is most
well-known for its central role in the proof of the 4 Color Theorem.  (For an
introduction to this technique, and a survey of results proved by it, see
\emph{A Guide to the Discharging Method}~\cite{discharging}, by the first author
and West.)  
We assume the theorem is false, and let $G$ be a counterexample with fewest
edges.  We assign to each vertex $v$ a charge
$d(v)-4$ and to each face $f$ a charge $\ell(f)-4$, where $d(v)$ and $\ell(f)$
denote the degree of $v$ and the length of $f$.  We denote these charges as
$\ch(v)$ and $\ch(f)$.
By Euler's formula, the sum of these initial charges (over all vertices and
faces) is $-8$, since $$\sum_{x\in V\cup F}\ch(x) = \sum_{v\in
V}d(v)-4+\sum_{f\in F}\ell(f)-4 = 2|E|-4|V|+2|E|-4|F| = -4 (2).
$$

Now we redistribute charge via the four discharging rules outlined below,
giving a final charge function $\ch^*$. Since $G$ is a minimal counterexample,
it must not contain any configurations that are reducible under either the 
\hyperref[basicreducibility]{Basic Reducibility Lemma} or the 
\hyperref[mainreducibility]{Main Reducibility Lemma}.
We use the absence of such configurations to show that each face and vertex
finishes with nonnegative final charge. This gives the following contradiction:
\[ -8 = \sum_{x \in V(G) \cup F(G)} \ch(x) = \sum_{x \in V(G) \cup F(G)} \ch^*(x) \ge 0. \]

Hence no such minimal counterexample $G$ can exist, so the Main Theorem is true.

\subsection{Discharging Rules}

\begin{figure}[h]
\centering
\begin{tikzpicture}
\coordinate (A) at (1.5, 0.0);
\coordinate (B) at (0.75, -1.29903);
\coordinate (C) at (-0.75, -1.29903);
\coordinate (D) at (-1.5, 0.0);
\coordinate (E) at (-0.75, 1.29903);
\coordinate (F) at (0.75, 1.29903);
\coordinate (G) at (-2.25, 1.29903);
\coordinate (H) at (2.25, 1.29903);
\coordinate (I) at (-3.75, 1.29903);
\coordinate (J) at (3.75, 1.29903);
\coordinate (K) at (-3.0, 0.0);
\coordinate (L) at (3.0, 0.0);
\coordinate (M) at (-2.25, -1.29903);
\coordinate (N) at (2.25, -1.29903);
\coordinate (O) at (0.75, -1.75);
\coordinate (P) at (-0.75, -1.75);

\draw (A)--(B)--(C)--(D)--(E)--(F)--(A) (E)--(G)--(D) (G)--(I)--(K)--(M)--(C) (F)--(H)--(A) (H)--(J)--(L)--(N)--(B) (B)--(O) (C)--(P);

\draw[fill] (A) circle [radius=\VertSize];
\draw[fill=white] (B) +(-\VertSize, \VertSize) rectangle +(\VertSize, -\VertSize) (B);
\draw[fill=white] (C) +(-\VertSize, \VertSize) rectangle +(\VertSize, -\VertSize) (C);
\draw[fill] (D) circle [radius=\VertSize];
\draw[fill] (E) circle [radius=\VertSize];
\draw[fill] (F) circle [radius=\VertSize];
\draw[fill] (G) circle [radius=\VertSize];
\draw[fill] (H) circle [radius=\VertSize];
\draw[fill] (I) circle [radius=\VertSize];
\draw[fill] (J) circle [radius=\VertSize];
\draw[fill] (K) circle [radius=\VertSize];
\draw[fill] (L) circle [radius=\VertSize];
\draw[fill] (M) circle [radius=\VertSize];
\draw[fill] (N) circle [radius=\VertSize];

\tikzset{arrows={-stealth}}
\draw (0.75, -0.43301)--(1.03125, -0.59538875)--(1.375, 0.0);
\draw (0.0, -0.86602)--(0.0, -1.1907775)--(0.65, -1.1907775);
\draw (0.0, -0.86602)--(0.0, -1.1907775)--(-0.65, -1.1907775);
\draw (-0.75, -0.43301)--(-1.03125, -0.59538875)--(-1.375, 0.0);
\draw (-0.75, 0.43301)--(-1.5, 0.86602);
\draw (0.0, 0.86602)--(0.0, 1.1907775)--(0.65, 1.1907775);
\draw (0.0, 0.86602)--(0.0, 1.1907775)--(-0.65, 1.1907775);
\draw (0.75, 0.43301)--(1.5, 0.86602);

\node at (0.65, -0.43301) {$ \frac{1}{3} $};
\node at (0.35, -0.86602) {$ \frac{1}{6} $};
\node at (-0.35, -0.86602) {$ \frac{1}{6} $};
\node at (-0.65, -0.43301) {$ \frac{1}{3} $};
\node at (-0.65, 0.43301) {$ \frac{1}{3} $};
\node at (0.35, 0.86602) {$ \frac{1}{6} $};
\node at (-0.35, 0.86602) {$ \frac{1}{6} $};
\node at (0.65, 0.43301) {$ \frac{1}{3} $};

\node at (-4.5, 0) {R1:};
\end{tikzpicture}
~~~~~
\begin{tikzpicture}
\coordinate (A) at (0, 0);
\coordinate (B) at (0.0, 1.5);
\coordinate (C) at (1.42658, 0.46352);
\coordinate (D) at (0.88167, -1.21352);
\coordinate (E) at (-0.88167, -1.21352);
\coordinate (F) at (-1.42658, 0.46352);

\draw (A)--(B) (A)--(C) (A)--(D) (A)--(E) (A)--(F);

\draw[fill] (A) circle [radius=\VertSize];
\draw[fill=white] (B) +(-\VertSize, \VertSize) rectangle +(\VertSize, -\VertSize) (B);
\draw[fill=white] (C) +(-\VertSize, \VertSize) rectangle +(\VertSize, -\VertSize) (C);
\draw[fill=white] (D) +(-\VertSize, \VertSize) rectangle +(\VertSize, -\VertSize) (D);
\draw[fill=white] (E) +(-\VertSize, \VertSize) rectangle +(\VertSize, -\VertSize) (E);
\draw[fill=white] (F) +(-\VertSize, \VertSize) rectangle +(\VertSize, -\VertSize) (F);

\node[right] at (B) {$ v_1 $};
\node[above] at (C) {$ v_2 $};
\node[right] at (D) {$ v_3 $};
\node[left] at (E) {$ v_4 $};
\node[above] at (F) {$ v_5 $};

\node at (-3, 0.5) {R2:};

\tikzset{arrows={-stealth}}
\draw (0.15, 0.2)--(1.25, 0.55);
\draw (0.25, -0.05)--(0.95, -1);
\draw (-0.25, -0.05)--(-0.95, -1);
\draw (-0.15, 0.2)--(-1.25, 0.55);

\node at (1, 0.75) {$ \frac{1}{4} $};
\node at (1, -0.6) {$ \frac{1}{4} $};
\node at (-1, -0.6) {$ \frac{1}{4} $};
\node at (-1, 0.75) {$ \frac{1}{4} $};
\end{tikzpicture}
~~~~~
\begin{tikzpicture}
\coordinate (A) at (0, 1);
\coordinate (B) at (-0.75, 0);
\coordinate (C) at (0.75, 0);

\draw (A)--(B)--(C)--(A);

\draw[fill] (A) circle [radius=\VertSize];
\draw[fill=white] (B) +(-\VertSize, \VertSize) rectangle +(\VertSize, -\VertSize) (B);
\draw[fill=white] (C) +(-\VertSize, \VertSize) rectangle +(\VertSize, -\VertSize) (C);

\node[right] at (A) {$ w $};
\node[left] at (B) {$ v $};
\node[right] at (C) {$ u $};

\node at (-1.5, 0.5) {R3:};

\tikzset{arrows={-stealth}}
\draw (-0.75, -0.25)--(0.75, -0.25)--(0, 0.75);
\node at (0, -0.5) {$ c $};
\end{tikzpicture}
~~~~~
\begin{tikzpicture}
\coordinate (A) at (0, 1);
\coordinate (X) at (0, 1.5);
\coordinate (Y) at (0, 2);
\coordinate (B) at (-0.75, 0);
\coordinate (C) at (0.75, 0);

\draw (A)--(B)--(C)--(A)--(X)--(Y);

\draw[fill] (A) circle [radius=\VertSize];
\draw[fill] (X) circle [radius=\VertSize];
\draw[fill=white] (Y) +(-\VertSize, \VertSize) rectangle +(\VertSize, -\VertSize) (Y);
\draw[fill=white] (B) +(-\VertSize, \VertSize) rectangle +(\VertSize, -\VertSize) (B);
\draw[fill=white] (C) +(-\VertSize, \VertSize) rectangle +(\VertSize, -\VertSize) (C);

\node[right] at (A) {$ w $};
\node[left] at (B) {$ v $};
\node[right] at (C) {$ u $};


\tikzset{arrows={-stealth}}
\draw (-0.75, -0.25)--(0.75, -0.25)--(0, 0.75);
\node at (0, -0.5) {$ c $};
\end{tikzpicture}
~~~~~
\begin{tikzpicture}
\coordinate (A) at (0, 0);
\coordinate (B) at (0.0, 1.0);
\coordinate (C) at (0.0, 2.0);
\coordinate (D) at (1.0, 0.0);
\coordinate (E) at (1.5, 0.86603);
\coordinate (F) at (1.5, -0.86603);
\coordinate (G) at (0.0, -1.5);
\coordinate (H) at (-1.5, 0.0);

\draw (A)--(B)--(C) (A)--(D)--(E) (D)--(F) (A)--(G) (A)--(H);

\draw[fill] (A) circle [radius=\VertSize];
\draw[fill] (B) circle [radius=\VertSize];
\draw[fill=white] (C) +(-\VertSize, \VertSize) rectangle +(\VertSize, -\VertSize) (C);
\draw[fill] (D) circle [radius=\VertSize];
\draw[fill=white] (E) +(-\VertSize, \VertSize) rectangle +(\VertSize, -\VertSize) (E);
\draw[fill=white] (F) +(-\VertSize, \VertSize) rectangle +(\VertSize, -\VertSize) (F);
\draw[fill=white] (G) +(-\VertSize, \VertSize) rectangle +(\VertSize, -\VertSize) (G);
\draw[fill=white] (H) +(-\VertSize, \VertSize) rectangle +(\VertSize, -\VertSize) (H);

\node[below left] at (A) {$ u $};
\node[right] at (B) {$ v_1 $};
\node[right] at (D) {$ v_2 $};
\node[right] at (G) {$ v_3 $};
\node[above] at (H) {$ v_4 $};

\node at (-3.3, 1) {R4:};

\tikzset{arrows={-stealth}}
\draw (-0.15, 0.15)--(-0.15, 0.85);
\draw (0.15, -0.15)--(0.85, -0.15);

\node at (-0.35, 0.5) {$ \frac{1}{4} $};
\node at (0.5, -0.5) {$ \frac{1}{4} $};
\end{tikzpicture}
\caption{
(R1) A 6-face gives charge $\frac13$ to each incident edge, and it is passed on
to either an incident 3-face or to one or both endpoints of the edge. (R2) If
$d(v_1) \ge 5$ and $d(v_i)\le 4$ for each $i\in\{2,3,4,5\}$, then the 5-vertex
splits its charge equally among $v_2$, $v_3$, $v_4$, and $v_5$.  (R3)  Here $u$
passes some or all of the charge it receives from $v$ on to $w$. (R4) If $u$
has charge $\frac12$ after R1--R3, then it splits this charge between its two
neighbors needing charge, $v_1$ and $v_2$.}
\label{fig:rules}
\end{figure}
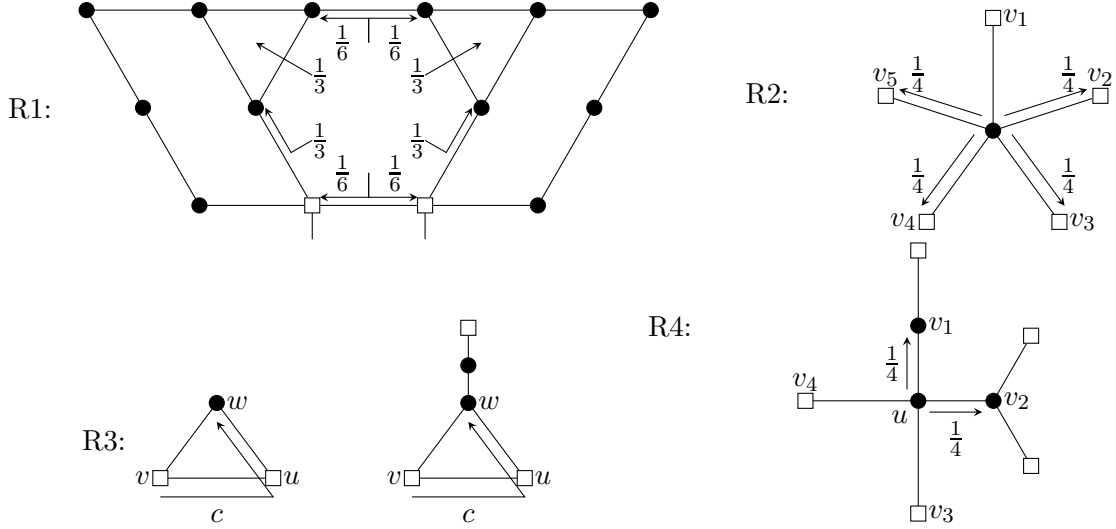

The following four discharging rules are applied to the elements of $ G $ successively, 
i.e., (R1) is applied everywhere that it is applicable, then (R2), then (R3),
and finally (R4).  Examples of these rules are illustrated in Figure~\ref{fig:rules}.
We write $k$-vertex (resp.~$k^+$, $k^-$) for a vertex of degree $k$ (resp.~at
least $k$, at most $k$).  We define $k$-faces analogously.

\begin{itemize}
\item[\textbf{R1}:] Each $ 6^+ $-face gives charge $ \frac{1}{3} $ to each incident edge. If such an edge $ e $ is incident to a 3-face $ f $, then $ e $ gives this charge to $ f $. Otherwise, $ e $ splits this charge evenly between any $ 3^- $-endpoints it has, or else splits it evenly between both endpoints if both have degree at least 4.
\footnote{Edges only ever act as a charge carrier between faces and other faces
or vertices. Outside of this phase, edges always have zero charge. Also, $G$
need not be 2-connected.  If a cutedge $e$ lies on a face $f$, then $f$ gives
$e$ charge $\frac23$.}

\item[\textbf{R2}:] 
Each $6^+$-vertex $v$ splits its initial charge evenly among its neighbors of
degree at most $d(v)$.
Each $5$-vertex with a $16^+$-neighbor splits its initial charge evenly among
its $4^-$-neighbors.
Each $5$-vertex $v$ with no $16^+$-neighbor splits is initial charge evenly among its
neighbors of the following types: 
3-vertices on triangular faces with $v$ and no $12^+$-neighbor, 2-vertices on
triangular faces with $v$, and other 2-vertices with no $( D-2)^+$-neighbor.

\item[\textbf{R3}:] Let $ u $ be a $ 4^+ $-vertex on a 3-face $ uvw $ and
suppose $ u $ receives some charge $ c $ during R2 from $ v $. If $ w $ is a
2-vertex, then $ u $ passes charge $ c $ on to $ w $. If instead $ w $ is a
3-vertex with a 2-neighbor whose other neighbor has degree less than $  D $, then $ u $ passes charge $ \min \{c, \frac{1}{2}\} $ on to $ w $.
\footnote{This rule rarely applies, and it 
can be largely ignored when seeking the high-level intuition behind the proof.}

\item[\textbf{R4}:] If a $ 3^+ $-vertex has positive charge after R1-R3, it splits this charge among its neighbors with negative charge, such that a 3-vertex gives charge at most $ \frac{4}{15} $ to another 3-vertex, and otherwise all charge splits evenly.
\end{itemize}

As stated above, we now show that $\ch^*(x) \ge 0$ for each vertex and face $x$.
It turns out that this is easy for everything except 3-vertices and
2-vertices, which require more detailed analysis.

\subsection{Faces and High-Degree Vertices}

All faces end with nonnegative final charge. Each $ 6^+ $-face $ f $ starts with
charge $ \ell(f) - 4 $ and gives away charge $ \frac{\ell(f)}{3} $. Thus $ f $
ends with $\ch^*(f) = \frac{2\ell(f)}{3} - 4 $, which is nonnegative since $
\ell(f) \ge 6 $. A 3-face cannot be adjacent to another 3-face since 4-cycles
are forbidden. Since $ G $ has no 4-cycles or 5-cycles, each 3-face $ f $ must
be adjacent to a $ 6^+ $-face on each of its edges. Each such $ 6^+ $-face
passes charge $ \frac{1}{3} $ to $ f $ via their common edge, so $\ch^*(f) = 3
- 4 + 3(\frac{1}{3}) = 0 $.

Each $ 4^+ $-vertex $ v $ starts out with nonnegative initial charge, and by the
design of the discharging rules never gives away more than its current charge,
so $\ch^*(v) \ge 0 $. Now we must verify that all 3-vertices and 2-vertices end
with nonnegative final charge as well, which will complete the proof.

\subsection{3-vertices}

First consider a 3-vertex $u$ that is not incident to any 3-faces. The three
faces meeting at $u$ must all be $ 6^+ $-faces, and thus each gives total
charge $\frac{2}{3}$ to two of the edges incident to $u$. Even when all of $
u $'s neighbors are $3^-$-vertices, $u$ receives at least half of this
charge, and hence end with $\ch^*(u) \ge 3 - 4 + 3(\frac{1}{3}) = 0$.

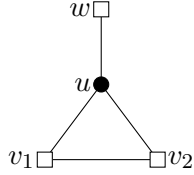
\begin{figure}[h]
\centering
\begin{tikzpicture}
\coordinate (A) at (0, 2);
\coordinate (B) at (0, 1);
\coordinate (C) at (-0.75, 0);
\coordinate (D) at (0.75, 0);

\draw (A)--(B)--(C)--(D)--(B);

\draw[fill=white] (A) +(-\VertSize, \VertSize) rectangle +(\VertSize, -\VertSize) (A);
\draw[fill] (B) circle [radius=\VertSize];
\draw[fill=white] (C) +(-\VertSize, \VertSize) rectangle +(\VertSize, -\VertSize) (C);
\draw[fill=white] (D) +(-\VertSize, \VertSize) rectangle +(\VertSize, -\VertSize) (D);

\node[left] at (A) {$ w $};
\node[left] at (B) {$ u $};
\node[left] at (C) {$ v_1 $};
\node[right] at (D) {$ v_2 $};
\end{tikzpicture}
\caption{The 3-vertex $ u $ on a 3-face under consideration.}
\label{fig:3vert_3face}
\end{figure}

Now consider a 3-vertex $ u $ on a 3-face $ uv_1v_2 $ whose third neighbor is $
w $, as shown in Figure~\ref{fig:3vert_3face}. Note that since $ v_1 $ and $
v_2 $ are adjacent, $ |N^2(u)| \le d(w) + d(v_1) + d(v_2) - 2 $. The two faces
incident to $ u $ other than the 3-face must be $ 6^+ $-faces, and hence 
give total charge $ \frac{2}{3} $ to the edge $ uw $. If $ d(w) \ge 4 $, then
all of this charge passes to $ u $, while if $ d(w) \le 3 $, then $ u $ 
receives charge $ \frac{1}{3} $ from this edge.

If $ d(v_i) = 2 $ for any $ i $, then $ v_i $ is reducible under the
\hyperref[basicreducibility]{Basic Reducibility Lemma}. Alternately, if $
d(v_i) \ge 12 $ for some $ i $, then $ u $ receives charge at least $
\frac{12 - 4}{12} = \frac{2}{3} $ from $ v_i $ via R2, and since $ uw $ sends $
u $ charge at least $ \frac{1}{3} $ via R1, $\ch^*(u) \ge 3 - 4 + \frac{2}{3} +
\frac{1}{3} = 0 $. Hence we can assume $ 3 \le d(v_i) \le 11 $ for $ i \in \{1,
2\} $. Also, if $ d(v_1) + d(v_2) \ge 17 $, then the
\hyperref[concavity]{Concavity Lemma} (with
$a=5$) implies that $u$
receives at least as much charge as when one $v_i$ is a $12^+$-vertex,
which, as just shown, ensures that $\ch^*(u) \ge 0$. Thus we assume
$d(v_1) + d(v_2) \le 16$. Now we consider what happens to $u$ based on the
degree of $w$.

\textit{Case $ d(w) \ge 6 $}: Here $u$ receives charge at least $\frac{6 -
4}{6} = \frac{1}{3}$ from $ w $, as well as charge $\frac{2}{3}$ from $uw$,
thus ends with $\ch^*(u) \ge 3 - 4 + \frac{2}{3} + \frac{1}{3} = 0$.

\textit{Case $ d(w) = 2 $}: Since $|N^2(u)| \le d(v_1) + d(v_2) \le 16$ and
$|N^2(w)| \le  D + 3$ and $ D \ge 14$, this configuration is
reducible under the \hyperref[mainreducibility]{Main Reducibility Lemma}.

\textit{Case $ d(w) \in \{3, 4, 5\} $}: We will show that $u$ receives charge
at least $\frac{1}{2}$ total from $w$ and the edge $uw$, and at least $
\frac{1}{4} $ from each of $ v_1 $ and $ v_2 $. This ensures that $\ch^*(u) \ge
3 - 4 + \frac{1}{2} + 2(\frac{1}{4}) = 0 $. First consider the charge from $ w
$ and $ uw $: if $ d(w) \ge 4 $, then as mentioned above, all $ \frac{2}{3} $
of the charge that passes through $ uw $ goes to $ u $, and $ \frac{2}{3} >
\frac{1}{2} $. Otherwise, if $ d(w) = 3 $, then $ u $ receives $ \frac{1}{3} $
from $ uw $, and so needs at least $ \frac{1}{6} $ more from $ w $ for this
total to reach $ \frac{1}{2} $.

Let $ x_1 $ and $ x_2 $ denote the neighbors of $ w $ other than $ u $. Since $
\{u, w\} $ is not reducible, the \hyperref[mainreducibility]{Main Reducibility
Lemma} implies that $ d(x_1) + d(x_2) \ge  D + 1 $. Now the
\hyperref[concavity]{Concavity Lemma}
implies that $w$ has at least as much charge to give to $u$ via R4 as 
when $d(x_1) =  D - 4$ and $d(x_2) = 5$ (and $x_2$ gives no charge to $u$).
If $w$ does not lie on a 3-face, then it receives charge $3(\frac{1}{3})$ from
its three incident edges via R1, making its charge nonnegative. Now the
additional charge of $ \frac{( D - 4) - 4}{ D - 4} $ from $x_1$ is 
split at most two ways. Since $ D \ge 10$, this ensures that $u$ gets an
additional charge of at least $\frac{1}{6}$ from $w$.

Suppose instead that $ w $ does lie on a 3-face. Now we know that $ d(x_2) \ge
3 $, since a 2-vertex on a 3-face with a 3-neighbor is reducible according to
the Basic Reducibility Lemma. Now if $ d(x_2) \ge 4 $, then $ x_2 $ always has
nonnegative charge and thus never needs to receive charge. If $ d(x_2) = 3 $,
then $ x_2 $ receives charge at least $ \frac{1}{3} $ from its incident
edge not on the 3-face, and at least $ \frac{2}{3} $ from $ x $ as long as $
d(x) \ge 12 $, meaning it does not need any charge from $ w $. Thus, whatever
the degree of $ x_2 $, vertex $w$ does not need to give it any charge via R4. 
Since $ D \ge 25$, this ensures that $w$ gets charge $\frac{1}{3}
+ \frac{5}{6}$ via R1 and R2, and thus gives charge $\frac{1}{6}$ to $u$ via
R4. Hence we have shown that $u$ always gets charge at least $\frac{1}{2}$ from
$w$ and the edge $uw$.
\smallskip

Now we show that $ u $ receives charge at least $ \frac{1}{4} $ from $ v_1 $
and, by symmetry, also from $ v_2 $. If $ d(v_1) \ge 6 $, then $ v_1 $ gives
charge at least $ \frac{1}{3} $ to $ u $ via R2, and $ \frac{1}{3} > \frac{1}{4}
$. Otherwise assume $ d(v_1) \le 5$.  

First suppose that $d(v_1)=5$.  If $v_1$ splits its charge between four or
fewer neighbors, then each receives charge at least $\frac14$, so we are done.
So assume instead that all five neighbors of $v_1$ should receive some of its
charge via R4.  We will show that $uv_1v_2$ is a reducible configuration.
By minimality, we can get a good ordering $\sigma'$ for $G-uv_2$.  
Let $S=\{v_2, u\}$. 
To extend $\sigma'$ to $G$, delete $S$ and append $v_2$, $u$; call this
ordering $\sigma$.  Clearly $\sigma$ is good for every vertex of $V(G)\setminus
S$.  Also, each vertex of $S$ has at most three neighbors in $G$ earlier in
$\sigma$.  Finally, each $x\in S$ has at most $ D+2$ neighbors in $G^2$
earlier in the ordering: 
$|N^2(v_2)\setminus\{u\}|\le d(v_1)+d(u)+( D-3)-3= D+2$ and
$|N^2(u)|\le 5+5+3$.  So assume $d(v_1)\in\{3,4\}$.

Recall that $ |N^2(u)| \le d(w) + d(v_1) + d(v_2) - 2 \le 19$. If $ \{u, v_1\}
$ is not reducible under the Main Reducibility Lemma, then $ |N^2(v_1)| \ge
 D + 4 $, i.e., $v_1$ has at least one high-degree neighbor $z$. 
Now $v_1$ has no excess charge to give to $u$ via R1, but will be able to give
the needed charge via R4. Note that by the same reasoning used above, since 
$\{u, v_2\}$ is not reducible under the \hyperref[mainreducibility]{Main
Reducibility Lemma}, $v_2$ must
also either be a $4^+$-vertex or have a high-degree neighbor. This means that $v_1$
never needs to
give charge to $v_2$ via R4, since $v_2$ only ever needs to receive charge if it
is a 3-vertex, and in such a case, it receives all the charge it needs from
its high-degree neighbor and incident edge off of the 3-face.

In the case that $d(v_1) = 3$, the neighbor $z$ of $v_1$ not on the 3-face
must have degree at least $ D - 8$.  Since $ D \ge 18$, this ensures
that $v_1$ gets charge at least $\frac{3}{5} + \frac{2}{3}$ from $z$ and
the edge $v_1 z$. Thus $v_1$ is able to pass charge at least 
$\frac{4}{15} > \frac{1}{4}$ to $u$.

\begin{figure}[h]
\centering
\begin{tikzpicture}
\coordinate (A) at (0, 2);
\coordinate (B) at (0, 1);
\coordinate (C) at (-0.75, 0);
\coordinate (D) at (0.75, 0);
\coordinate (E) at (-2.25, 0);
\coordinate (F) at (-1.5, 1);
\coordinate (G) at (-1.5, 1.5);
\coordinate (H) at (-1.5, 2);

\draw (A)--(B)--(C)--(E)--(F)--(C)--(D)--(B) (F)--(G)--(H);

\draw[fill=white] (A) +(-\VertSize, \VertSize) rectangle +(\VertSize, -\VertSize) (A);
\draw[fill] (B) circle [radius=\VertSize];
\draw[fill] (C) circle [radius=\VertSize];
\draw[fill=white] (D) +(-\VertSize, \VertSize) rectangle +(\VertSize, -\VertSize) (D);
\draw[fill=white] (E) +(-\VertSize, \VertSize) rectangle +(\VertSize, -\VertSize) (E);
\draw[fill] (F) circle [radius=\VertSize];
\draw[fill] (G) circle [radius=\VertSize];
\draw[fill=white] (H) +(-\VertSize, \VertSize) rectangle +(\VertSize, -\VertSize) (E);

\node[left] at (A) {$ w $};
\node[left] at (B) {$ u $};
\node[below] at (C) {$ v_1 $};
\node[right] at (D) {$ v_2 $};
\node[left] at (E) {$ z $};
\node[left] at (F) {$ t $};
\node[left] at (G) {$ s $};
\end{tikzpicture}
\caption{This configuration, where R3 would apply, is reducible by the Main
Reducibility Lemma.}
\label{fig:no_R3}
\end{figure}
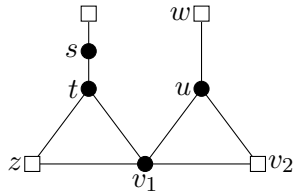

If instead $ d(v_1) = 4 $, then $ v_1 $ split any excess charge it
receives at most two ways via R4 (since neither $ z $ nor $ v_2 $ needs
charge). Let $t$ be the neighbor of $v_1$ other than $u, v_2$, and $z$, and
note that $v_1$ only sends charge to $t$ via R4 if $d(t) < 4$. By the
\hyperref[concavity]{Concavity Lemma}, $v_1$ receives no less charge than when
$d(z) =  D - 5
$, $d(t) = 3$, and $d(v_2) = 5$ (but it doesn't give charge to $v_1$).
If $v_1zt$ is not a 3-face, then $v_1$
receives charge at least $\frac{( D - 5) - 4}{ D - 5} +
\frac{1}{3}$ from $z$ and the edge $v_1z$. Since $ D \ge 10$, 
vertex $v_1$ get charge at least $\frac{8}{15}$, so it passes at least
$\frac{4}{15} > \frac{1}{4}$ to $u$ via R4.

If instead $ v_1zt $ is a 3-face, then we note that $ t $ cannot be a 2-vertex,
since this would be reducible. Also, $t$ cannot be a 3-vertex with a
2-neighbor $ s $, where the other neighbor of $ s $ has degree less than $
 D $, because this also would be reducible under the
\hyperref[mainreducibility]{Main Reducibility Lemma} (using the vertex sequence
$S = \{t, s, u\}$), as shown in
Figure~\ref{fig:no_R3}. Since these are the only times when R3 can apply, we
conclude that this rule is not used here. Hence $ v_1 $ gets charge at least $
\frac{( D - 5) - 4}{ D - 5}$ from $z$, which it can then send at
least half of to $u$. As long as $ D \ge 13$, this means $v_1$ sends
at least $\frac{1}{4}$ to $u$ as desired.

\subsection{2-vertices}

\begin{figure}[h]
\centering
\begin{tikzpicture}
\coordinate (A) at (0, 0);
\coordinate (B) at (1, 1);
\coordinate (C) at (2, 0);

\draw (A)--(B)--(C)--(A);

\tikzset{arrows={-stealth}}
\draw (0, 0.2)--(0.8, 1);
\draw (0.4, 0.2)--(1.6, 0.2)--(1, 0.8);
\draw (2, 0.2)--(1.2, 1);

\draw[fill=white] (A) +(-\VertSize, \VertSize) rectangle +(\VertSize, -\VertSize) (A);
\draw[fill] (B) circle [radius=\VertSize];
\draw[fill=white] (C) +(-\VertSize, \VertSize) rectangle +(\VertSize, -\VertSize) (C);

\node[below] at (A) {$ v_1 $};
\node[above] at (B) {$ u $};
\node[below] at (C) {$ v_2 $};
\end{tikzpicture}
\caption{A 2-vertex on a 3-face receives charge via R2 and R3.}
\label{fig:2vert_R3}
\end{figure}
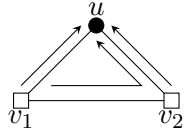

\textbf{2-vertex on a 3-face:} First consider a 2-vertex $ u $ on a 3-face $
uv_1v_2 $, as depicted in Figure~\ref{fig:2vert_R3}. By the Basic Reducibility
Lemma, this is reducible unless $ d(v_1) + d(v_2) \ge  D + 5 $. By the
\hyperref[concavity]{Concavity Lemma}, we know that $ u $ receives at least as much charge 
as if $d(v_1) =  D$ and $d(v_2) = 5$. Now $u$ receives charge at least $
\frac{ D - 4}{ D} + \frac{1}{4} $ via R2. However, $ v_2 $ also
receives charge $ \frac{ D - 4}{ D} $ from $ v_1 $ via R2, and the
conditions are met for R3, so $ v_2 $ passes this charge along to $ u $.
Hence in total $ u $ receives charge at least $ 2(\frac{ D - 4}{ D}) +
\frac{1}{4} $. Since $ D \ge 32$, $u$ ends with $\ch^*(u) \ge 2 - 4
+ 2(\frac{32 - 4}{32}) + \frac{1}{4} = 0$.

\textbf{2-vertex with one high-degree neighbor:} Now we assume that the
2-vertex $u$, with neighbors $v_1$ and $v_2$, does not lie on a 3-face.
Note that if $d(v_i) = 2$ for some $i \in \{1, 2\}$, then $ \{u, v_i\} $ is
reducible under the \hyperref[mainreducibility]{Main Reducibility Lemma}. Hence
we assume that $d(v_1) \ge 3$ and $d(v_2) \ge 3$. 

Suppose $d(v_1) \ge  D - 2$; now $u$ receives charge $\frac{2}{3}$ through
the edge $u v_1$ via R1 and $\frac{( D - 2) - 4}{ D - 2}$ from $v_1$
via R2. If $d(v_2) \ge 4$, then $u$ also gets $\frac{2}{3}$ through the edge $u
v_2$ via R1, and so ends with final charge at least $ 2 - 4 + 2(\frac{2}{3}) +
\frac{( D - 2) - 4}{ D - 2} $, which is nonnegative since $ D \ge 14$.

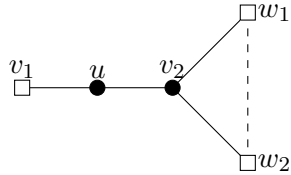
\begin{figure}[h]
\centering
\begin{tikzpicture}
\coordinate (A) at (-1, 0);
\coordinate (B) at (0, 0);
\coordinate (C) at (1, 0);
\coordinate (D) at (2, 1);
\coordinate (E) at (2, -1);

\draw (A)--(B)--(C)--(D) (C)--(E);

\draw[dashed] (D)--(E);

\draw[fill=white] (A) +(-\VertSize, \VertSize) rectangle +(\VertSize, -\VertSize) (A);
\draw[fill] (B) circle [radius=\VertSize];
\draw[fill] (C) circle [radius=\VertSize];
\draw[fill=white] (D) +(-\VertSize, \VertSize) rectangle +(\VertSize, -\VertSize) (D);
\draw[fill=white] (E) +(-\VertSize, \VertSize) rectangle +(\VertSize, -\VertSize) (E);

\node[above] at (A) {$ v_1 $};
\node[above] at (B) {$ u $};
\node[above] at (C) {$ v_2 $};
\node[right] at (D) {$ w_1 $};
\node[right] at (E) {$ w_2 $};
\end{tikzpicture}
\caption{A 2-vertex $ u $ with a neighbor $ v_1 $ such that $ d(v_1) \ge  D - 2 $.}
\label{2vert_big_nbr}
\end{figure}

So assume $d(v_2) = 3$, and denote the other neighbors of $v_2$ by $w_1$
and $w_2$, as pictured in Figure~\ref{2vert_big_nbr}. Note that $v_2$ and $u$
each receive charge $\frac{1}{3}$ from the edge $u v_2$ via R1. Now $\{u,
v_2\}$ is reducible under the \hyperref[mainreducibility]{Main
Reducibility Lemma} unless $ |N^2(v_2)| \ge  D + 3 $. First, suppose that $v_2$
lies on a 3-face, which implies $ d(w_1) + d(w_2) \ge  D + 3 $. By the
\hyperref[concavity]{Concavity Lemma}, $ v_2 $ receives at least as much charge
as if $ d(w_1) =  D - 1 $ and $ d(w_2) = 4 $. Hence after R2, $ v_2 $ has
charge at least $ 3 - 4 + \frac{1}{3} + \frac{( D - 1) - 4}{ D - 1} $.
Since $  D \ge 26 $, this ensures that $ v_2 $ has charge at least $ -1 +
\frac{1}{3} + \frac{21}{25} > \frac{1}{6} $ after R2, which it passes to $ u $
via R4. (Note that $ w_2 $ does not receive charge from $ v_2 $ via R4: since $
v_2w_1w_2 $ is a 3-face, $ d(w_2) > 2 $. Further, if $ d(w_2) = 3 $, then $ w_2
$ receives enough charge from $ w_1 $ and its incident edge off of the
3-face.) Hence $\ch^*(u) \ge 2 - 4 + \frac{2}{3} + \frac{1}{3} + \frac{(26 - 2)
- 4}{26 - 2} + \frac{1}{6} = 0 $.

So suppose instead that $ v_2 $ does not lie on a 3-face. Now $ |N^2(v_2)| \ge
 D + 3 $, implying that $ d(w_1) + d(w_2) \ge  D + 1 $. Again using the
\hyperref[concavity]{Concavity Lemma}, we can assume that $ d(w_1) \ge  D -
4$. Now $v_2$ gets charge at least $\frac{1}{3}$ from each of the edges $
u v_2 $ and $ v_2 w_2 $, and $ \frac{2}{3} $ from the edge $ v_2 w_1 $ via R1,
which already puts its total charge at $ 3 - 4 + \frac{4}{3} = \frac{1}{3} $.
Now $ v_2 $ splits this charge at most two ways (giving to $u$ and
possibly $w_2$) via R4. Since $ v_2 $ has charge at least $ \frac{1}{3} $
after R1, it gives charge at least $ \frac{1}{6} $ to $ u $ via R4. As
shown above, since $  D \ge 26 $ this ensures that $\ch^*(u) \ge 0 $, as
desired.

Hereafter we assume that $ d(v_1) \le  D - 3 $ and $ d(v_2) \le  D - 3
$. We show that $ u $ must receive total charge at least 1 from edge $ u v_1 $
and vertex $ v_1 $; by symmetry the same is true of edge $ u v_2 $ and vertex $
v_2 $. This ensures that $ u $ ends with final charge at least $ 2 - 4 + 1 + 1
= 0 $, as desired. If $ d(v_1) \ge 6 $, then $ u $ gets charge $ \frac{2}{3} $
from $ u v_1 $ via R1 and charge $ \frac{d(v_1) - 4}{d(v_1)} \ge \frac{6 -
4}{6} = \frac{1}{3} $ from $ v_1 $ via R2. This gives $ u $ the charge of 1
from $ v_1 $'s side as needed, so henceforth we assume $ d(v_1) \le 5 $.
\bigskip

\textbf{2-vertex with a 3-neighbor:} Suppose $ d(v_1) = 3 $, and denote the
other neighbors of $v_1$ by $w_1$ and $w_2$, with $d(w_1) \ge d(w_2)$.
Now $u$ receives charge $ \frac{1}{3} $ from the edge $ u v_1 $ via R1,
meaning it needs to get $ \frac{2}{3} $ from $ v_1 $ via R4. First suppose that
$ v_1 $ does not lie on a 3-face.
Since $ d(v_2) \le  D - 3 $, we apply the \hyperref[mainreducibility]{Main
Reducibility Lemma} with $ S = \{v_1, u\} $, unless $ d(w_1) + d(w_2) \ge
 D + 2 $. Likewise, if $ d(w_2) = 2 $, then we simply take $ S = \{v_1, w_2, u\} $.

Hence we assume $ d(w_2) \ge 3 $. If $ d(w_2) \ge 4 $, then $ v_1 $ receives
charge $ \frac{2}{3} $ from both of the edges $ v_1 w_1 $ and $ v_1 w_2 $,
along with $ \frac{1}{3} $ from the edge $ u v_1 $ via R1. This means that
after R1 alone, $ v_1 $ has charge $ 3 - 4 + \frac{1}{3} + 2(\frac{2}{3})
= \frac{2}{3} $, which it can then send to $ u $ via R4 as needed. So instead
suppose that $ d(w_2) = 3 $, which implies $ d(w_1) \ge  D - 1 $. Now $ v_1 $
gets charge at least $ \frac{4}{3} $ via R1 ($ \frac{1}{3} $ each from edges $
u v_1 $ and $ v_1 w_2 $, and $ \frac{2}{3} $ from edge $ v_1 w_1 $) and $
\frac{( D - 1) - 4}{ D - 1} $ from $ w_1 $ via R2. Since $  D \ge 11
$, this ensures that $ v_1 $ has charge at least $ 3 - 4 + \frac{4}{3} +
\frac{(11 - 1) - 4}{11 - 1} = \frac{14}{15} $ after R2. Since $ v_1 $ gives no
more charge than $ \frac{4}{15} $ to $ w_2 $ via R4, it can give at least $
\frac{10}{15} = \frac{2}{3} $ to $ u $ via R4 as needed. So $ u $ gets charge
at least 1 from $ v_1 $ and $ uv_1 $.

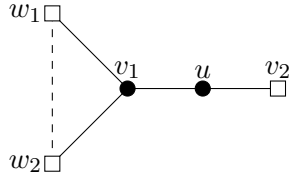
\begin{figure}[h]
\centering
\begin{tikzpicture}
\coordinate (A) at (-2, 1);
\coordinate (B) at (-2, -1);
\coordinate (C) at (-1, 0);
\coordinate (D) at (0, 0);
\coordinate (E) at (1, 0);

\draw (A)--(C)--(D)--(E) (B)--(C);

\draw[dashed] (A)--(B);

\draw[fill=white] (A) +(-\VertSize, \VertSize) rectangle +(\VertSize, -\VertSize) (A);
\draw[fill=white] (B) +(-\VertSize, \VertSize) rectangle +(\VertSize, -\VertSize) (B);
\draw[fill] (C) circle [radius=\VertSize];
\draw[fill] (D) circle [radius=\VertSize];
\draw[fill=white] (E) +(-\VertSize, \VertSize) rectangle +(\VertSize, -\VertSize) (E);

\node[left] at (A) {$ w_1 $};
\node[left] at (B) {$ w_2 $};
\node[above] at (C) {$ v_1 $};
\node[above] at (D) {$ u $};
\node[above] at (E) {$ v_2 $};
\end{tikzpicture}
\caption{A 2-vertex $ u $ with a 3-neighbor $ v_1 $.}
\label{2vert_3nbr}
\end{figure}

Now suppose instead that $ v_1 $ does lie on a 3-face. If we cannot apply the
\hyperref[mainreducibility]{Main Reducibility Lemma} with $ S = \{v_1, u\} $,
then $ d(w_1) + d(w_2) \ge  D + 4 $. By the \hyperref[concavity]{Concavity
Lemma}, $ v_1 $ receives at least as much charge as if $ d(w_1) =  D $
and $ d(w_2) = 4 $. Thus $ v_1 $ receives charge $ \frac{1}{3} $ from edge $ u
v_1 $ via R1, and further receives charge at least $ \frac{ D - 4}{ D}
$ from $ w_1 $ via R2. Additionally, $ w_2 $ receives at least $ \frac{ D -
4}{ D} $ from $ w_1 $ via R2, and the criteria are met for R3; since $
 D \ge 8 $, this means $ w_2 $ passes charge $ \frac{1}{2} $ to $ v_1 $.
Hence after R3, $ v_1 $ has charge at least $ 3 - 4 + \frac{1}{3} + \frac{1}{2}
+ \frac{ D - 4}{ D} $. Since $  D \ge 24 $, this means $ v_1 $ 
has charge at least $ -\frac{1}{6} + (\frac{24 - 4}{24}) = \frac{2}{3} $ that
it can pass to $ u $ via R4, as needed.
\bigskip

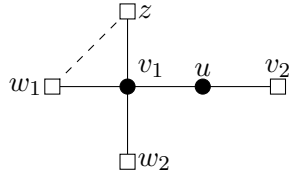
\begin{figure}[h]
\centering
\begin{tikzpicture}
\coordinate (A) at (-1, 1);
\coordinate (B) at (-2, 0);
\coordinate (C) at (-1, 0);
\coordinate (D) at (-1, -1);
\coordinate (E) at (0, 0);
\coordinate (F) at (1, 0);

\draw (A)--(C)--(E)--(F) (B)--(C)--(D);

\draw[dashed] (A)--(B);

\draw[fill=white] (A) +(-\VertSize, \VertSize) rectangle +(\VertSize, -\VertSize) (A);
\draw[fill=white] (B) +(-\VertSize, \VertSize) rectangle +(\VertSize, -\VertSize) (B);
\draw[fill] (C) circle [radius=\VertSize];
\draw[fill=white] (D) +(-\VertSize, \VertSize) rectangle +(\VertSize, -\VertSize) (D);
\draw[fill] (E) circle [radius=\VertSize];
\draw[fill=white] (F) +(-\VertSize, \VertSize) rectangle +(\VertSize, -\VertSize) (F);

\node[right] at (A) {$ z $};
\node[left] at (B) {$ w_1 $};
\node[above right] at (C) {$ v_1 $};
\node[right] at (D) {$ w_2 $};
\node[above] at (E) {$ u $};
\node[above] at (F) {$ v_2 $};
\end{tikzpicture}
\caption{A 2-vertex $ u $ with a 4-neighbor $ v_1 $, where $ v_1 $ has a high-degree neighbor $ z $.}
\label{2vert_4nbr}
\end{figure}

\textbf{2-vertex with a 4-neighbor:} Now suppose $ d(v_1) = 4 $. In this case, $
u $ receives charge $ \frac{2}{3} $ from edge $ u v_1 $ via R1, and hence only
needs to get charge $ \frac{1}{3} $ more from $ v_1 $ via R4. We can apply the
\hyperref[mainreducibility]{Main Reducibility Lemma} with $ S = \{v_1, u\} $
unless $ |N^2(v_1)| \ge  D + 4 $, which means the degree sum of the
neighbors of $ v_1 $ other than $ u $ is at least $  D + 2 $. The least
charge that passes from $v_1$ to $u$ via R4 occurs when $v_1$ has as many
$3^-$-neighbors as possible, so we assume that $ v_1 $ has two $3^-$-neighbors
$ w_1 $ and $ w_2 $ and one high-degree neighbor $ z $, as shown in
Figure~\ref{2vert_4nbr}.

By the \hyperref[concavity]{Concavity Lemma}, $v_1$ receives at least as much
charge via R2 as if $d(z) =  D - 8 $ and $d(w_1) = d(w_2) = 5$ (but neither
$w_1$ nor $w_2$ gives charge to $v_1$). If $v_1$ and $z$ do not lie on a
common 3-face, then $ v_1 $ receives charge $\frac{1}{3}$ from edge $ v_1 z $
via R1. Since $  D \ge 20$, $ v_1 $
receives charge at least $ \frac{(20 - 8) - 4}{20 - 8} = \frac{2}{3} $ from $ z
$ via R2, giving $ v_1 $ a total charge of at least 1 after R2. Since $ v_1 $
splits its charge at most three ways, it passes charge at least $
\frac{1}{3} $ to $ u $ via R4, as needed.

Instead, assume $ v_1 z w_1 $ is a 3-face. By the
\hyperref[basicreducibility]{Basic Reducibility Lemma}, we know $ w_1 $ cannot
be a 2-vertex, so instead assume $ d(w_1) \ge 3 $. 
First, suppose $d(w_1)=3$, and let $x$ be the third
neighbor of $ w_1 $ besides $ v_1 $ and $ z $. Now $ w_1 $ receives charge at
least $ \frac{1}{3} $ from edge $ w_1x $ via R1 and, since $  D \ge 20 $,
receives charge at least $ \frac{(20 - 8) - 4}{20 - 8} = \frac{2}{3} $ from $ z
$ via R2. Hence $ w_1 $ has nonnegative charge after R2, and thus does not need
charge from $ v_1 $ via R4, meaning $ v_1 $ only splits its charge at
most two ways.  Similarly, if $d(w_1)\ge 4$, then $w_1$ does not need charge
from $v_1$ via R4.  Thus, in every case, $v_1$ splits its charge after R3 at most
two ways. 

Now $ v_1 $ also receives charge at least $ \frac{2}{3} $ from $ z $ via R1. If
$ d(x) = 2 $ and the other neighbor of $ x $ has degree less than $  D $,
then the sequence $ S = \{w_1, x, u\} $ is reducible under the
\hyperref[mainreducibility]{Main Reducibility Lemma}. If instead $ d(x) \ge 3
$, or $ d(x) = 2 $ and the other neighbor of $ x $ has degree $  D $, then
the conditions for R3 are not met, which means $ v_1 $ keeps its charge from $
z $ until R4. Splitting at most two ways, $ v_1 $ can give charge at least $
\frac{1}{3} $ to $ u $ via R4, which is all $ u $ still needs.
\bigskip

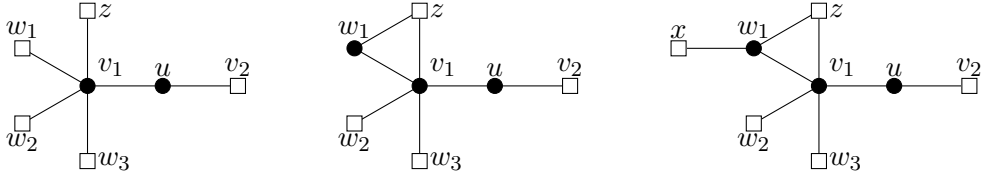
\begin{figure}[h]
\centering
\begin{tikzpicture}
\coordinate (Z) at (-1, 1);
\coordinate (W1) at (-1.86603, 0.5);
\coordinate (V1) at (-1, 0);
\coordinate (W2) at (-1.86603, -0.5);
\coordinate (W3) at (-1, -1);
\coordinate (U) at (0, 0);
\coordinate (V2) at (1, 0);

\draw (Z)--(V1)--(U)--(V2) (W1)--(V1)--(W3) (V1)--(W2);

\draw[fill=white] (Z) +(-\VertSize, \VertSize) rectangle +(\VertSize, -\VertSize) (Z);
\draw[fill=white] (W1) +(-\VertSize, \VertSize) rectangle +(\VertSize, -\VertSize) (W1);
\draw[fill] (V1) circle [radius=\VertSize];
\draw[fill=white] (W2) +(-\VertSize, \VertSize) rectangle +(\VertSize, -\VertSize) (W2);
\draw[fill=white] (W3) +(-\VertSize, \VertSize) rectangle +(\VertSize, -\VertSize) (W3);
\draw[fill] (U) circle [radius=\VertSize];
\draw[fill=white] (V2) +(-\VertSize, \VertSize) rectangle +(\VertSize, -\VertSize) (V2);

\node[right] at (Z) {$ z $};
\node[above] at (W1) {$ w_1 $};
\node[above right] at (V1) {$ v_1 $};
\node[below] at (W2) {$ w_2 $};
\node[right] at (W3) {$ w_3 $};
\node[above] at (U) {$ u $};
\node[above] at (V2) {$ v_2 $};
\end{tikzpicture}
~~~~~
\begin{tikzpicture}
\draw (Z)--(V1)--(U)--(V2) (W1)--(V1)--(W3) (V1)--(W2) (Z)--(W1);

\draw[fill=white] (Z) +(-\VertSize, \VertSize) rectangle +(\VertSize, -\VertSize) (Z);
\draw[fill] (W1) circle [radius=\VertSize];
\draw[fill] (V1) circle [radius=\VertSize];
\draw[fill=white] (W2) +(-\VertSize, \VertSize) rectangle +(\VertSize, -\VertSize) (W2);
\draw[fill=white] (W3) +(-\VertSize, \VertSize) rectangle +(\VertSize, -\VertSize) (W3);
\draw[fill] (U) circle [radius=\VertSize];
\draw[fill=white] (V2) +(-\VertSize, \VertSize) rectangle +(\VertSize, -\VertSize) (V2);

\node[right] at (Z) {$ z $};
\node[above] at (W1) {$ w_1 $};
\node[above right] at (V1) {$ v_1 $};
\node[below] at (W2) {$ w_2 $};
\node[right] at (W3) {$ w_3 $};
\node[above] at (U) {$ u $};
\node[above] at (V2) {$ v_2 $};
\end{tikzpicture}
~~~~~
\begin{tikzpicture}
\coordinate (X) at (-2.86603, 0.5);

\draw (Z)--(V1)--(U)--(V2) (W1)--(V1)--(W3) (V1)--(W2) (Z)--(W1)--(X);

\draw[fill=white] (Z) +(-\VertSize, \VertSize) rectangle +(\VertSize, -\VertSize) (Z);
\draw[fill] (W1) circle [radius=\VertSize];
\draw[fill=white] (X) +(-\VertSize, \VertSize) rectangle +(\VertSize, -\VertSize) (X);
\draw[fill] (V1) circle [radius=\VertSize];
\draw[fill=white] (W2) +(-\VertSize, \VertSize) rectangle +(\VertSize, -\VertSize) (W2);
\draw[fill=white] (W3) +(-\VertSize, \VertSize) rectangle +(\VertSize, -\VertSize) (W3);
\draw[fill] (U) circle [radius=\VertSize];
\draw[fill=white] (V2) +(-\VertSize, \VertSize) rectangle +(\VertSize, -\VertSize) (V2);

\node[right] at (Z) {$ z $};
\node[above] at (W1) {$ w_1 $};
\node[above] at (X) {$ x $};
\node[above right] at (V1) {$ v_1 $};
\node[below] at (W2) {$ w_2 $};
\node[right] at (W3) {$ w_3 $};
\node[above] at (U) {$ u $};
\node[above] at (V2) {$ v_2 $};
\end{tikzpicture}
\caption{Cases where a 2-vertex $ u $ has a 5-neighbor $ v_1 $.}
\label{2vert_5nbr}
\end{figure}

\textbf{2-vertex with a 5-neighbor:} Finally, suppose $d(v_1) = 5$, as shown in
Figure~\ref{2vert_5nbr}. Similar to above, $u$ receives charge
$\frac{2}{3}$ from edge $uv_1$ via R1.  Now we must consider whether or not
$v_1$ has a $16^+$-neighbor.  First, suppose that it does.

Since $v_1$ has a $16^+$-neighbor,
it splits its initial charge of $5 - 4 = 1$ at most four ways, so it
passes charge at least $\frac{1}{4}$ to $u$ via R2. Thus in order for $u$ to
receive charge at least 1 from $v_1$ and the edge $uv_1$, it only needs to get
charge $\frac{1}{12}$ more from $v_1$ via R4.

Let $z$ denote the highest-degree neighbor of $v_1$, and denote its other
neighbors by $w_1$, $w_2$, and $w_3$. If $v_1$ and $z$ are not together on
a 3-face, then $v_1$ receives charge $\frac{1}{3}$ from edge $v_1z$
via R1, and does not lose this charge prior to R4. Thus in R4, $v_1$ has
charge at least $\frac{1}{3}$ which it splits at most four ways, meaning it
sends charge at least $\frac{1}{12}$ to $u$, as needed. So instead assume
that $v_1zw_1$ is a 3-face.  Now since $|N^2(v_1)| \ge  D + 4$, we
have $d(z) + d(w_1) + d(w_2) + d(w_3) \ge  D + 4$; by the
\hyperref[concavity]{Concavity Lemma}, $v_1$ receives at least as much charge
via R1 and R2 as if $d(z) =  D - 10$ and $d(w_1)=4$ and  $d(w_2) = d(w_3) =
5$ (but neither sends charge to $v_1$ via R2).

Suppose $d(w_1) = 2$. This configuration is not immediately reducible under
either the \hyperref[basicreducibility]{Basic Reducibility Lemma} or the
\hyperref[mainreducibility]{Main Reducibility Lemma}, but is in
fact reducible using a hybrid of the two approaches. If we delete vertex
$w_1$ as in the \hyperref[basicreducibility]{Basic Reducibility Lemma}, we get
a good ordering $\sigma'$ for 
$G-w_1$.  To extend this ordering to $G$, we delete $u$ and append $w_1, u$.
The key point is that now $u$ is not an earlier neighbor of $w_1$ in $G^2$,
so the number of earlier neighbors for $w_1$ in $G^2$ is at most
$d(z)+d(v_1)-2-1\le  D+5-3= D+2$.  Also, recall that we
are assuming $d(v_2)\le  D-3$, so $|N^2(u)|\le d(v_1)+d(v_2)\le
5+( D-3)= D+2$.
Hence, this configuration is reducible.

Now assume $d(w_1) \ge 3$. If $d(w_1) \ge 4$ then whatever charge $v_1$
gets from $ z $ via R2 it keeps until R4. Since $d(z)\ge 6$, this means that
$v_1$ receives charge at least $\frac{6 - 4}6 = \frac{1}{3} $ in
R2, and splits it at most three ways in R4, so it gives $u$ charge at
least $\frac{1}{9} > \frac{1}{12}$.  Instead suppose $d(w_1) = 3$, and let
$x$ be the other neighbor of $w_1$. If the criteria for R3 are not met
(i.e.\ $d(x) \ge 3$ or $d(x) = 2$ and the other neighbor of $x$ has degree
$ D$), then $v_1$ keeps any charge it receives from $z$ via R2 until R4.
Thus, as before, $v_1$ still gets charge at least $\frac{1}{3}$ since $d(z)\ge
6$, and splitting at most four ways gives charge $\frac{1}{12}$ to $u$
via R4, as needed.

Suppose instead that $d(x) = 2$ and the other neighbor of $x$ has degree at most $
 D - 1$.  Now $v_1$ passes some charge that it gets from $z$ via R2 to
$w_1$ via R3. Since $d(z)\ge 16$, $v_1$ receives charge at least $\frac{16 -
4}{16} = \frac{3}{4}$ from $z$ via R2. Now $v_1$
gives charge $\frac{1}{2}$ to $w_1$ via R3, leaving it with charge $\frac{3}{4}
- \frac{1}{2} = \frac{1}{4}$. Since $w_1$ gets charge at least $\frac{1}{3}$
from the edge $w_1x$ via R1, $\frac{3}{4}$ from $z$ via
R2, and $\frac{1}{2}$ from $ v_1 $ via R3, it has nonnegative charge, and
thus needs no charge from $ v_1 $ via R4. Hence $v_1$ splits its remaining
$\frac{1}{4}$ charge at most three ways, meaning it gives charge at least $
\frac{1}{12}$ to $u$ via R4 as needed.
\bigskip

Now suppose instead that $v_1$ has no $16^+$-neighbor.
Since $u$ receives charge $\frac23$ from edge $uv_1$, we must show that in this
case $u$ still receives charge at least $\frac13$ from vertex $v_1$.  By R2,
$v_1$ splits its charge of 1 among neighbors of the following types: 
3-vertices on triangular faces with $v_1$ and no $12^+$-neighbor, 2-vertices on
triangular faces with $v_1$, and other 2-vertices with no $( D-2)^+$-neighbor.
If $v_1$ has at most three neighbors of these types,
then clearly $v_1$ gives charge at least $\frac13$ to $u$, and we are done.
So, suppose instead that $v_1$ has at least four neighbors of these types.  
In particular, this implies that $v_1$ has at most one $4^+$-neighbor and no
$16^+$-neighbor.
We will show that $G$ contains a reducible configuration.

Note that $v_1$ can be incident to at most two triangular faces.  
We will show that $v_1$ gives charge via R2 to at most neighbor not on a triangular face 
and at most one neighbor on each of at most two incident triangular faces.  
Thus, $v_1$ gives charge to at most 3 neibhbors by R2.

Suppose that $v_1$ has two 2-neighbors, say $u_1$ and $u_2$, such that
each $u_i$ has no $( D-2)^+$-neighbor.  Form $G'$ from $G$ by deleting $u_1$
and $u_2$.  By minimality, $G'$ has a good vertex ordering $\sigma'$. 
To reach a good vertex ordering $\sigma$ for $G$, delete $v_1$ from $\sigma'$,
then append $v_1$, $u_1$, $u_2$.  Now $v_1$ has at most three earlier
neighbors in $\sigma$ and at most $15+(2)3+(2)1=22$ earlier neighbors in $G^2$. 
Also, each $u_i$ has at most two earlier neighbors in $G$ and at most
$( D-3)+5$ earlier neighbors in $G^2$.

Now we must verify that on each incident triangular face $v_1$ has at most one
neighbor that receives charge.  If $v_1$ has two such neighbors on a common 3-face
and one is a 2-neighbor, say $u_2$, then the configuration is reducible by the
\hyperref[basicreducibility]{Basic Reducibility Lemma}, since $|N^2(u_2)|\le
5+3$.  So suppose that $v_1$ has
two 3-neighbors, $u_2$ and $u_3$, on a common 3-face and they both receive
charge from $v_1$.  Form $G'$ from $G$ by deleting edge $u_2u_3$.  By
minimality, $G'$ has a good vertex ordering $\sigma'$.  To get a good vertex
ordering $\sigma$ for $G$, delete $u_2$ and $u_3$ from $\sigma'$, then append
$u_2$ and $u_3$.  Clearly, each $u_i$ has at most 3 earlier neighbors in the
ordering.  Also, $v_1$ gives charge to $u_2$ only when $u_2$ has no
$12^+$-neighbor.  Thus, $|N^2(u_2)|\le 5+3+11$; similarly for $u_3$.  Thus, the
resulting vertex ordering $\sigma$ is good for $G$.
$\hfill \qed$
\bigskip

To conclude the paper, we remark that this vertex ordering guaranteed by the
\hyperref[main]{Main Theoerm} can be constructed 
recursively in linear time.  The basic idea is to find some
reducible configuration in amortized constant time.  We assume a data structure
that stores for each vertex: its degree, a doubly-linked adjacency list in
clockwise order, and for
each neighbor a pointer to that neighbor.  Note that to handle each reducible
configuration, we either delete a vertex of low degree or we delete an edge with
both endpoints of low degree.  Thus, we can preprocess $G$ in linear time to
find all such reducible configurations, storing them in some generic ``bag''
(for example a stack or a queue).  Now at each step, we remove some reducible
configuration from the bag, recurse on the appropriate smaller graph, and add to
the bag any newly created reducible configurations.  (The proof of the
\hyperref[main]{Main Theorem} guarantees that the bag will never be empty.)
The first author and
Kim give a lengthier explanation of these ideas in Section 6
of~\cite{cranston08}.

\section*{Acknowledgment}
Thanks to Landon Rabern for his careful reading of the manuscript, which caught
some inaccuracies.

\section*{Appendix}
\label{appendix}
In this section, we first collect a few standard graph theory definitions.  We
conclude with a construction of Dvo\v{r}\'{a}k et al.~\cite{dvorak} of planar
graphs $G$ of girth 6 and maximum degree $\Delta$ such that $\chi(G^2)\ge
\Delta+2$ (for each $\Delta \ge 2$).

The \emph{girth} of a graph is the length of its shortest cycle.  
The degree $d(v)$ of a vertex $v$ is its number of incident
edges.  The maximum degree in $G$ is denoted $\Delta$.  The set of vertices
within distance 2 of a vertex $v$ is denoted $N^2(v)$.
We write $k$-vertex (resp.~$k^+$, $k^-$) for a vertex of degree $k$ (resp.~at
least $k$, at most $k$).  We define $k$-faces analogously.

A \emph{coloring} of a graph $G$ assigns to each vertex a color (typically denoted by a
positive integer).  
A coloring $f$ is \emph{proper} if the endpoints $u$ and
$v$ of each edge $uv$ get distinct colors, i.e., $f(u)\ne f(v)$.  A graph is
\emph{$k$-colorable} if it has a proper coloring
with at most $k$ colors.  The \emph{chromatic number} $\chi(G)$ of a graph $G$
is the least $k$ such that $G$ is $k$-colorable.
A \emph{list assignment} $L$ assigns to each vertex $v$ a set of allowable
colors $L(v)$.  An \emph{$L$-coloring} is a proper coloring $f$ such that
$f(v)\in L(v)$ for every vertex $v$.  A graph $G$ is \emph{$k$-choosable} if it
is $L$-colorable whenever $|L(v)|=k$ for every $v\in V(G)$.  The \emph{list
chromatic number $\chil(G)$ of $G$} (or \emph{choice number} of $G$) is the
least $k$ such that $G$ is $k$-choosable.

The game of {$k$-paintability} (or online list $k$-coloring) is played by two
players, \emph{Lister} and \emph{Painter}.  In each round $i$, Lister presents
to Painter some nonempty list (set) of uncolored vertices.  Painter chooses
(paints) some subset of them to receive color $i$.  If Lister lists some
particular vertex $k$ times and Painter never paints it, then Lister wins. 
Otherwise Painter wins.  The \emph{paint number} $\chi_p(G)$ is the least $k$
such that $G$ is $k$-paintable.

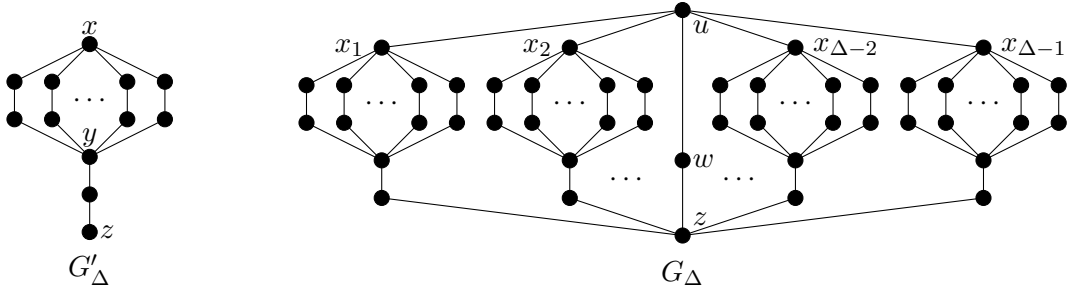
\begin{figure}[h]
\centering
\begin{tikzpicture}
\coordinate (A) at (0, 1.5);
\coordinate (B) at (-1, 1);
\coordinate (C) at (-0.5, 1);
\coordinate (D) at (0.5, 1);
\coordinate (E) at (1, 1);
\coordinate (F) at (-1, 0.5);
\coordinate (G) at (-0.5, 0.5);
\coordinate (H) at (0.5, 0.5);
\coordinate (I) at (1, 0.5);
\coordinate (J) at (0, 0);
\coordinate (K) at (0, -0.5);
\coordinate (Z) at (0, -1);

\draw (Z)--(K)--(J)--(F)--(B)--(A)--(C)--(G)--(J)--(H)--(D)--(A)--(E)--(I)--(J);

\draw[fill] (A) circle [radius = \VertSize];
\draw[fill] (B) circle [radius = \VertSize];
\draw[fill] (C) circle [radius = \VertSize];
\draw[fill] (D) circle [radius = \VertSize];
\draw[fill] (E) circle [radius = \VertSize];
\draw[fill] (F) circle [radius = \VertSize];
\draw[fill] (G) circle [radius = \VertSize];
\draw[fill] (H) circle [radius = \VertSize];
\draw[fill] (I) circle [radius = \VertSize];
\draw[fill] (J) circle [radius = \VertSize];
\draw[fill] (K) circle [radius = \VertSize];
\draw[fill] (Z) circle [radius = \VertSize];

\node at (0, 0.75) {$ \cdots $};

\node[above] at (A) {$ x $};
\node[above] at (J) {$ y $};
\node[right] at (Z) {$ z $};

\node at (0, -1.5) {$ G'_\Delta $};
\end{tikzpicture}
~~~~~~~~~~
\begin{tikzpicture}
\coordinate (U) at (0, 2);
\coordinate (V) at (0, -1);
\coordinate (W) at (0, 0);

\coordinate (A1) at (-4, 1.5);
\coordinate (B1) at (-5, 1);
\coordinate (C1) at (-4.5, 1);
\coordinate (D1) at (-3.5, 1);
\coordinate (E1) at (-3, 1);
\coordinate (F1) at (-5, 0.5);
\coordinate (G1) at (-4.5, 0.5);
\coordinate (H1) at (-3.5, 0.5);
\coordinate (I1) at (-3, 0.5);
\coordinate (J1) at (-4, 0);
\coordinate (K1) at (-4, -0.5);

\coordinate (A2) at (-1.5, 1.5);
\coordinate (B2) at (-2.5, 1);
\coordinate (C2) at (-2, 1);
\coordinate (D2) at (-1, 1);
\coordinate (E2) at (-0.5, 1);
\coordinate (F2) at (-2.5, 0.5);
\coordinate (G2) at (-2, 0.5);
\coordinate (H2) at (-1, 0.5);
\coordinate (I2) at (-0.5, 0.5);
\coordinate (J2) at (-1.5, 0);
\coordinate (K2) at (-1.5, -0.5);

\coordinate (A3) at (1.5, 1.5);
\coordinate (B3) at (0.5, 1);
\coordinate (C3) at (1, 1);
\coordinate (D3) at (2, 1);
\coordinate (E3) at (2.5, 1);
\coordinate (F3) at (0.5, 0.5);
\coordinate (G3) at (1, 0.5);
\coordinate (H3) at (2, 0.5);
\coordinate (I3) at (2.5, 0.5);
\coordinate (J3) at (1.5, 0);
\coordinate (K3) at (1.5, -0.5);

\coordinate (A4) at (4, 1.5);
\coordinate (B4) at (3, 1);
\coordinate (C4) at (3.5, 1);
\coordinate (D4) at (4.5, 1);
\coordinate (E4) at (5, 1);
\coordinate (F4) at (3, 0.5);
\coordinate (G4) at (3.5, 0.5);
\coordinate (H4) at (4.5, 0.5);
\coordinate (I4) at (5, 0.5);
\coordinate (J4) at (4, 0);
\coordinate (K4) at (4, -0.5);

\draw (V)--(K1)--(J1)--(F1)--(B1)--(A1)--(C1)--(G1)--(J1)--(H1)--(D1)--(A1)--(E1)--(I1)--(J1) (U)--(A1);
\draw (V)--(K2)--(J2)--(F2)--(B2)--(A2)--(C2)--(G2)--(J2)--(H2)--(D2)--(A2)--(E2)--(I2)--(J2) (U)--(A2);
\draw (V)--(K3)--(J3)--(F3)--(B3)--(A3)--(C3)--(G3)--(J3)--(H3)--(D3)--(A3)--(E3)--(I3)--(J3) (U)--(A3);
\draw (V)--(K4)--(J4)--(F4)--(B4)--(A4)--(C4)--(G4)--(J4)--(H4)--(D4)--(A4)--(E4)--(I4)--(J4) (U)--(A4);

\draw (U)--(W)--(V);

\draw[fill] (U) circle [radius = \VertSize];
\draw[fill] (V) circle [radius = \VertSize];
\draw[fill] (W) circle [radius = \VertSize];

\draw[fill] (A1) circle [radius = \VertSize];
\draw[fill] (B1) circle [radius = \VertSize];
\draw[fill] (C1) circle [radius = \VertSize];
\draw[fill] (D1) circle [radius = \VertSize];
\draw[fill] (E1) circle [radius = \VertSize];
\draw[fill] (F1) circle [radius = \VertSize];
\draw[fill] (G1) circle [radius = \VertSize];
\draw[fill] (H1) circle [radius = \VertSize];
\draw[fill] (I1) circle [radius = \VertSize];
\draw[fill] (J1) circle [radius = \VertSize];
\draw[fill] (K1) circle [radius = \VertSize];

\draw[fill] (A2) circle [radius = \VertSize];
\draw[fill] (B2) circle [radius = \VertSize];
\draw[fill] (C2) circle [radius = \VertSize];
\draw[fill] (D2) circle [radius = \VertSize];
\draw[fill] (E2) circle [radius = \VertSize];
\draw[fill] (F2) circle [radius = \VertSize];
\draw[fill] (G2) circle [radius = \VertSize];
\draw[fill] (H2) circle [radius = \VertSize];
\draw[fill] (I2) circle [radius = \VertSize];
\draw[fill] (J2) circle [radius = \VertSize];
\draw[fill] (K2) circle [radius = \VertSize];

\draw[fill] (A3) circle [radius = \VertSize];
\draw[fill] (B3) circle [radius = \VertSize];
\draw[fill] (C3) circle [radius = \VertSize];
\draw[fill] (D3) circle [radius = \VertSize];
\draw[fill] (E3) circle [radius = \VertSize];
\draw[fill] (F3) circle [radius = \VertSize];
\draw[fill] (G3) circle [radius = \VertSize];
\draw[fill] (H3) circle [radius = \VertSize];
\draw[fill] (I3) circle [radius = \VertSize];
\draw[fill] (J3) circle [radius = \VertSize];
\draw[fill] (K3) circle [radius = \VertSize];

\draw[fill] (A4) circle [radius = \VertSize];
\draw[fill] (B4) circle [radius = \VertSize];
\draw[fill] (C4) circle [radius = \VertSize];
\draw[fill] (D4) circle [radius = \VertSize];
\draw[fill] (E4) circle [radius = \VertSize];
\draw[fill] (F4) circle [radius = \VertSize];
\draw[fill] (G4) circle [radius = \VertSize];
\draw[fill] (H4) circle [radius = \VertSize];
\draw[fill] (I4) circle [radius = \VertSize];
\draw[fill] (J4) circle [radius = \VertSize];
\draw[fill] (K4) circle [radius = \VertSize];

\node at (-4, 0.75) {$ \cdots $};
\node at (-1.5, 0.75) {$ \cdots $};
\node at (1.5, 0.75) {$ \cdots $};
\node at (4, 0.75) {$ \cdots $};
\node at (-0.75, -0.25) {$ \cdots $};
\node at (0.75, -0.25) {$ \cdots $};

\node[left] at (-4.1, 1.5) {$ x_1 $};
\node[left] at (-1.6, 1.5) {$ x_2 $};
\node[right] at (1.6, 1.5) {$ x_{\Delta - 2} $};
\node[right] at (4.1, 1.5) {$ x_{\Delta - 1} $};

\node[below right] at (U) {$ u $};
\node[above right] at (V) {$ z $};
\node[right] at (W) {$ w $};

\node at (0, -1.5) {$ G_\Delta $};
\end{tikzpicture}
\caption{In any $ (\Delta + 1) $-coloring of the square of $ G'_\Delta $, the $ (\Delta - 1) $-vertex $ x $ and the 1-vertex $ z $ cannot receive the same color. Because of this, no $ (\Delta + 1) $-coloring of the square of $ G_\Delta $ is possible, hence $ \chi(G_\Delta^2) \ge \Delta + 2 $.}
\label{lowerbounds}
\end{figure}

Now we present a construction of planar graphs $G_\Delta$ with maximum degree
$\Delta$ and girth 6 such that $\chi(G^2_\Delta)\ge \Delta+2$.  The first such
construction appeared in Borodin et al.~\cite{girthseven}.  The construction
we present is due to 
Dvo\v{r}\'{a}k et al~\cite{dvorak}.  We like it because we find it simpler, and
the graphs it produces have fewer vertices.

The key to the construction is a gadget $G_{\Delta}'$, show on the left in
Figure~\ref{lowerbounds}.  It consists of two vertices $x$ and $y$ joined by
$\Delta-1$ paths of length 3, as well as another path of length 2 incident to
vertex $y$; call the other endpoint of this 2-path $z$.  The key observation is that
in any coloring of $(G_{\Delta}')^2$ with $\Delta+1$ colors, vertices $x$ and $z$
must receive distinct colors.  The reason is that $y$ and all of its neighbors
must receive the $\Delta+1$ distinct colors.  So $z$ must receive the same color
as some neighbor $t$ of $y$ other than its common neighbor with $z$.  This
neighbor $t$
will be distance 2 from $x$, so it cannot recieve the same color as $x$.
To form $G_\Delta$, we take $\Delta-1$ copies of the gadget, identifying vertex
$z$ in all of them.  Further, we add a new vertex $u$ adjacent to $x$ in each
gadget, and we add a new vertex $w$ adjacent to $u$ and $z$.  Now the vertex set
$\{u,w,z,x_1,\ldots, x_{\Delta-1}\}$ has size $\Delta+2$ and in a coloring
of $G^2$ each pair of its vertices must receive distinct colors.  Thus,
$\chi(G^2_\Delta)\ge \Delta+2$.



\begin{thebibliography}{99}
\bibitem{agnarsson}
G. Agnarsson and M. Halld\'{o}rsson,
\textit{Coloring powers of planar graphs},
SIAM J. Discrete Math. \textbf{16} (2003), no. 4, 651--662.




\bibitem{BBGH1}
O. V. Borodin, H. J. Broersma, A. Glebov, and J. van den Heuvel,
{Stars and bunches in planar graphs, Part I: Triangulations},
\emph{CDAM Research Report Series} {\bf 2002--04} (2002). Original Russian version.
\emph{Diskretn. Anal. Issled. Oper. Ser. 1} {\bf 8} (2001), no. 2, 15--39,
available at:
\url{http://www.cdam.lse.ac.uk/Reports/Abstracts/cdam-2002-04.html}.

\bibitem{BBGH2}
O.V. Borodin, H. J. Broersma, A. Glebov, and J. van den Heuvel,
{Stars and bunches in planar graphs, Part II: General planar graphs and
colourings},
\emph{CDAM Research Report Series} {\bf 2002--05} (2002). Original Russian version.
\emph{Diskretn. Anal. Issled. Oper. Ser. 1} {\bf 8} (2001), no. 4, 9--33,
available at:
\url{http://www.cdam.lse.ac.uk/Reports/Abstracts/cdam-2002-05.html}.

\bibitem{girthseven} O.V. Borodin, A.N. Glebov, A.O. Ivanova, T.K. Neustroeva, and V.A. Tashkinov, \textit{A sufficient condition for a planar graph to be 2-distant $ (\Delta + 1) $-colorable}, Sib. Elektron. Math. Izv. \textbf{1} (2004), 129-141 [Russian].

\bibitem{girthsix2} O.V. Borodin and A.O. Ivanova, \textit{2-distance $ (\Delta + 2) $-coloring of planar graphs with girth six and $ \Delta \ge 18 $}, Discrete Math. \textbf{309} (2009), 6496-6502.



\bibitem{cranston08}
D.W. Cranston and S.-J. Kim, 
\textit{List-coloring the square of a subcubic graph}, 
J. Graph Theory \textbf{57} (2008), no. 1, 65--87. 

\bibitem{discharging} D.W. Cranston and D.B. West, \textit{A guide to the
discharging method}, preprint, available at: \url{http://arxiv.org/abs/1306.4434}.

\bibitem{dvorak} Z. Dvo\v{r}\'{a}k, D. Kr\'{a}l', P. Nejedl\'{y}, and R.
\v{S}krekovski, \textit{Coloring squares of planar graphs with girth six},
European J. Combin. \textbf{29} (2008), no. 4, 838-849.


\bibitem{havet} F. Havet, J. van den Heuvel, C. McDiarmid, and B. Reed,
\textit{List colouring squares of planar graphs}, preprint, available at:
\url{http://arxiv.org/abs/0807.3233}.

\bibitem{vandenheuvel} J. van den Heuvel and S. McGuinness, \textit{Coloring
the square of a planar graph}, J. Graph Theory \textbf{42} (2003), 110--124.



\bibitem{wanglih} W.-F. Wang, and K.-W. Lih, 
\textit{Labeling planar graphs with conditions on girth and distance two}, 
SIAM J. Discrete Math. 17 (2003), no. 2, 264--275.

\bibitem{wegner} G. Wegner, \textit{Graphs with given diameter and a colouring problem}, preprint, University of Dortmund (1977).

\bibitem{zhu-etal}
H.-Y. Zhu, X.-Z. Lu, C.-Q. Wang, and M. Chen, 
\textit{Labeling planar graphs without 4,5-cycles with a condition on distance
two.} SIAM J. Discrete Math. 26 (2012), no. 1, 52--64. 

\end{thebibliography}
\end{document}